\magnification=1200
\catcode`\@=11
\hsize=125mm  \vsize=187mm
\hoffset=4mm  \voffset=10mm
\pretolerance=500  \tolerance=1000  \brokenpenalty=5000

\catcode`\;=\active
\def;{\relax\ifhmode\ifdim\lastskip>\z@
\unskip\fi\kern.2em\fi\string;}

\catcode`\:=\active
\def:{\relax\ifhmode\ifdim\lastskip>\z@\unskip\fi
\penalty\@M\ \fi\string:}

\catcode`\!=\active
\def!{\relax\ifhmode\ifdim\lastskip>\z@
\unskip\fi\kern.2em\fi\string!}

\catcode`\?=\active
\def?{\relax\ifhmode\ifdim\lastskip>\z@
\unskip\fi\kern.2em\fi\string?}

\def\^#1{\if#1i{\accent"5E\i}\else{\accent"5E#1}\fi}
\def\"#1{\if#1i{\accent"7F\i}\else{\accent"7F#1}\fi}

\frenchspacing

\def\hexnumber@#1{\ifnum#1<10 \number#1\else
 \ifnum#1=10 A\else\ifnum#1=11 B\else\ifnum#1=12 C\else
 \ifnum#1=13 D\else\ifnum#1=14 E\else\ifnum#1=15 F\fi\fi\fi\fi\fi\fi\fi}
\def\mathhexbox@#1#2#3{\text{$\m@th\mathchar"#1#2#3$}}
\def\text{\relaxnext@\ifmmode\let\next\text@\else\let\next\text@@\fi\next}
\def\text@@#1{\leavevmode\hbox{#1}}
\font\tenmsa=msam10
\font\sevenmsa=msam7
\font\fivemsa=msam5
\font\tenmsb=msbm10
\font\sevenmsb=msbm7
\font\fivemsb=msbm5
\font\teneuf=eufm10
\font\seveneuf=eufm7
\font\fiveeuf=eufm5
\newfam\msafam
\newfam\msbfam
\newfam\euffam
\textfont\msafam=\tenmsa
\scriptfont\msafam=\sevenmsa
\scriptscriptfont\msafam=\fivemsa
\textfont\msbfam=\tenmsb
\scriptfont\msbfam=\sevenmsb
\scriptscriptfont\msbfam=\fivemsb
\textfont\euffam=\teneuf
\scriptfont\euffam=\seveneuf
\scriptscriptfont\euffam=\fiveeuf
\def\msa@{\hexnumber@\msafam}
\def\msb@{\hexnumber@\msbfam}
\def\Bbb{\relaxnext@\ifmmode\let\next\Bbb@\else
 \def\next{\errmessage{Utilisez \string\Bbb\space seulement en mode
math}}\fi\next}
\def\Bbb@#1{{\Bbb@@{#1}}}
\def\Bbb@@#1{\fam\msbfam#1}
\def\goth{\relaxnext@\ifmmode\let\next\goth@\else
 \def\next{\errmessage{Utilisez \string\goth\space seulement en mode
math}}\fi\next}
\def\goth@#1{{\goth@@{#1}}}
\def\goth@@#1{\fam\euffam#1}
\def\relaxnext@{\let\next\relax}

\catcode`\@=12

\font\eightrm=cmr8
\newif\ifpagetitre  \pagetitretrue
\newtoks\hautpagetitre \hautpagetitre={\hfil}
\newtoks\baspagetitre  \baspagetitre={\hfil}
\newtoks\auteurcourant  \auteurcourant={\hfil}
\newtoks\titrecourant  \titrecourant={\hfil}
\newtoks\hautpagegauche  \newtoks\hautpagedroite
\hautpagegauche={\llap{
\oldstyle\folio\quad}\eightrm\the\auteurcourant\hfill}
\hautpagedroite={\hfill\eightrm\the\titrecourant\quad\oldstyle\folio}
\newtoks\baspagegauche  \baspagegauche={\hfil} \newtoks\baspagedroite
\baspagedroite={\hfil} \headline={\ifpagetitre\the\hautpagetitre
\else\ifodd\pageno\the\hautpagedroite \else\the\hautpagegauche\fi\fi}
\footline={\ifpagetitre\the\baspagetitre
\global\pagetitrefalse
\else\ifodd\pageno\the\baspagedroite
\else\the\baspagegauche\fi\fi}

\def\ent{\mathop{\rm E\kern 0pt}\nolimits}

\def\lin{\mathop{\cal L\kern 0pt}\nolimits}

\def\build#1_#2{\mathrel{\mathop{\kern 0pt#1}\limits_{#2}}}
\def\tend_#1{\build\hbox to 1cm{\rightarrowfill}_{#1}}

\def\maj#1#2:{#1{\eightrm #2 : }}

\def\diagram#1{\def\normalbaselines{\baselineskip=0pt
\lineskip=10pt\lineskiplimit=1pt}  \matrix{#1}}

\def\hfl#1#2{\smash{\mathop{\hbox to 12mm{\rightarrowfill}}
\limits^{\scriptstyle#1}_{\scriptstyle#2}}}

\def\vfl#1#2{\llap{$\scriptstyle #1$}\left\downarrow
\vbox to 6mm{}\right.\rlap{$\scriptstyle #2$}}

\def\buildrel#1\over#2{\mathrel{
\mathop{\kern 0pt#2}\limits^{#1}}}

\def\build#1_#2^#3{\mathrel{
\mathop{\kern 0pt#1}\limits_{#2}^{#3}}}

 2
 4

\newfam\bffam \textfont\bffam=\tenbf \scriptfont\bffam=\sevenbf
\scriptscriptfont\bffam=\fivebf
\def\bf{\fam\bffam\tenbf}

\def\boxit#1#2{\setbox1=\hbox{\kern#1{#2}\kern#1}%
\dimen1=\ht1 \advance\dimen1 by #1 \dimen2=\dp1 \advance\dimen2 by #1
\setbox1=\hbox{\vrule height\dimen1 depth\dimen2\box1\vrule}%
\setbox1=\vbox{\hrule\box1\hrule}%
\advance\dimen1 by .4pt \ht1=\dimen1
\advance\dimen2 by .4pt \dp1=\dimen2 \box1\relax}

\overfullrule=0pt

\def\un{\mathop{\cup}\limits}
\def\CAP{\mathop{\cap}\limits}
\def\INF{\mathop{\inf}\limits}
\def\SUP{\mathop{\sup}\limits}
\def\fleche{\longrightarrow}

\def\''#1#2#3{\hbox{$\eta_{#3,#1+#2}-\eta_{#3,#1}-\eta_{#3,#2}$}\ }
\def\X#1#2#3{\hbox{$\xi_{#1,#2,#3}$}}
\def\XX#1#2#3{\X{#1}{#2}{#3}-\X{#1}{#3}{#2}+\X{#3}{#1}{#2}}

\def\h'#1#2#3{\XX{#1}{#2}{#3}=\ \''{#1}{#2}{#3}\ }     

\def\bs{\bigskip}

\def\ms{\medskip}

\def\nat{\hbox{{\rm N}\kern-1em\hbox{{\rm I}}}}

\def\C{\hbox{{\rm C}\kern-.55em\hbox{{\rm I}}}}

\def\Z{\ \hbox{{\rm Z}\kern-.4em\hbox{\rm Z}}}

\def\dis{\displaystyle}

\def\boxit1#1#2{\vbox{\hrule\hbox{\vrule
\vbox spread#1{\vfil\hbox spread#1{\hfil#2\hfil}\vfil}%
\vrule}\hrule}}
\def\boxit#1#2{\setbox1=\hbox{\kern#1{#2}\kern#1}%
\dimen1=\ht1 \advance\dimen1 by #1 \dimen2=\dp1 \advance\dimen2 by #1
\setbox1=\hbox{\vrule height\dimen1 depth\dimen2\box1\vrule}%
\setbox1=\vbox{\hrule\box1\hrule}%
\advance\dimen1 by .4pt \ht1=\dimen1
\advance\dimen2 by .4pt \dp1=\dimen2 \box1\relax}

\def\bs{\bigskip}

\def\ms{\medskip}

\def\trait{\centerline{\hbox to 3cm{\ \hrulefill\ }}}

\def\build#1_#2^#3{\mathrel{\mathop{\kern 0pt#1}\limits_{#2}^{#3}}}

\def\hfl#1#2{\smash{\mathop{\hbox to 12mm{\rightarrowfill}}
\limits^{\scriptstyle#1}_{\scriptstyle#2}}}

\def\vfl#1#2{\llap{$\scriptstyle #1$}\left\downarrow
\vbox to 6mm{}\right.\rlap{$\scriptstyle #2$}}

\def\hlf#1#2{\smash{\mathop{\hbox to 12mm{\leftarrowfill}}
\limits^{\scriptstyle#1}_{\scriptstyle#2}}}

\def\hfll#1#2{\smash{\mathop{\hbox to 12mm{\rightarrowfill}}
\limits^{\scriptstyle\hbox{$#1$}}_{\scriptstyle\hbox{$#2$}}}}

\def\hllf#1#2{\smash{\mathop{\hbox to 12mm{\leftarrowfill}}
\limits^{\scriptstyle\hbox{$#1$}}_{\scriptstyle\hbox{$#2$}}}}

\def\diagram#1{\def\normalbaselines{\baselineskip=0pt
\lineskip=10pt\lineskiplimit=1pt}  \matrix{#1}}

\pagetitrefalse

\hsize=125mm  \vsize=187mm
\hoffset=4mm  \voffset=10mm
\baselineskip=15pt
\pretolerance=500  \tolerance=1000  \brokenpenalty=5000
\def\un{\mathop{\cup}\limits}
\def\CAP{\mathop{\cap}\limits}
\def\CUP{\mathop{\cup}\limits}
\def\fleche{\longrightarrow}
\def\im{\mathop{\rm Im}\limits}
\def\ext{\mathop{\rm Ext}\limits}
\def\int{\mathop{\rm Int}\limits}

\bs\bs\bs
{\centerline {{\bf TH\'EORIE DE LA  MESURE DANS LES LIEUX R\'EGULIERS}}}
\bs
{\centerline{ {\bf ou   }  }}

\bs
{\centerline { {\bf LES INTERSECTIONS CACH\'EES  DANS LE PARADOXE DE BANACH-TARSKI}}} 

\bs\bs

{\centerline{ Olivier Leroy}}
\vfill

\eject
{\centerline{PREAMBULE}}
\bs
La disparition brutale d'Olivier Leroy en avril 1996 a 
priv\'e la communaut\'e math\'ematique d'un esprit brillant\footnote{(*)}{voir par
exemple : A. Grothendieck, R\'ecoltes et Semailles, Tome II, note 96, pages 409 \`a 413.}
 mais bien m\'econnu.
\ms
Le texte qui suit est la transcription fid\`ele d'un manuscrit qu'il avait
r\'edig\'e en 1995 \`a la suite de deux expos\'es donn\'es au s\'eminaire AGATA
sur le probl\`eme de la mesure dans les lieux.
D'autres travaux sur le m\^eme sujet ont \'et\'e  laiss\'es en chantier par 
O.Leroy, 
en particulier une amorce de th\'eorie de l'int\'egration, mais 
 nous ne pr\'esentons
ici, pour l'instant, que la partie consacr\'ee \`a la mesure.   
\ms
Nous avons pris l'initiative de compl\'eter le texte original par : 
une introduction, les notes en bas de page, des ajouts
entre doubles crochets ([[\dots]]), ainsi qu'une annexe 
(sous-espaces et sous-lieux d'un espace topologique). 
\bs
\ms
\bs\bs
\hfill Montpellier le  20 mai  1998
\bs
\hfill Jean Malgoire, Christine Voisin
 \vfil
\eject
{\centerline {\bf INTRODUCTION}}
\bs 
Il est bien connu que l'axiome du choix implique l'existence de parties non
mesurables pour la mesure de Lebesgue sur $\Bbb R$ ainsi que l'existence de 
d\'ecompositions ``paradoxales" de la boule unit\'e de $\Bbb R^3$ (Banach-Tarski).
Ceci est g\'en\'eralement interpr\'et\'e comme le prix \`a payer pour les
nombreux services  rendus par cet axiome \footnote{(*)}{par exemple, 
comme nous l'avait
fait remarquer Olivier Leroy, la moyennabili-\break t\'e des groupes ab\'eliens  : i.e.
l'existence de mesures additives (ou densit\'e) invariantes par translations,
de masse totale 1, d\'efinies sur l'ensemble de {\sl toutes\/} les parties.}.
\ms
La th\'eorie propos\'ee par Olivier Leroy montre que l'on peut avoir 
simultan\'ement l'axiome du choix et ``tout est mesurable"\ ! Elle se place
dans le cadre des {\sl lieux\/} (terme fran\c cais choisi par l'auteur comme 
traduction de l'anglais {\sl locales}) qui sont des
 cas particuliers  de topos au
sens de Grothendieck : un lieu est simplement un ensemble ordonn\'e qui a les 
propri\'et\'es formelles de l'ensemble ordonn\'e des ouverts d'un espace 
topologique. Les lieux ont d\'ej\`a fait l'objet de nombreuses \'etudes (cf 
\footnote{(**)}{Sheaves in Geometry and Logic. S.Mac Lane \& I.Moerdijk.
Springer 92.} et son abondante bibliographie).
      
\ms
Un des aspects remarquables de cette th\'eorie (en particulier pour ceux, 
nombreux, qui ont une m\'efiance syst\'ematique pour les g\'en\'eralisations...)
est qu'elle s'applique de mani\`ere pertinente  aux espaces topologiques habituels
dans lesquels elle fait appara\^itre
 des ``sous-espaces non classiques" 
(des sous-lieux);
avec pour cons\'equence que l'intersection (aux sens des sous-lieux) de
 sous-espaces ordinaires n'est plus forc\'ement un sous-espace ordinaire.
 \footnote{(***)}{on a par exemple un sous-lieu de $\Bbb R$ 
(appel\'e {\sl lieu g\'en\'erique de $\Bbb R$ })
 qui est l'intersection des ouverts denses et qui est encore dense
(bien que sans point!).}  
\ms
Le r\'esultat le plus frappant est sans doute que le prolongement naturel
de la mesure ext\'erieure de Lebesgue (sur $[0,1]$ par exemple) \`a {\sl tous\/}
les sous-lieux de $[0,1]$  est une {\sl mesure $\sigma$-additive
ext\'erieurement r\'eguli\`ere\/} !!!
\ms
 Les partitions (ensemblistes) ``paradoxales"
 qui donnaient des parties\break  non mesurables au sens de Lebesgue ne sont
 plus des partitions au sens des lieux : il y a des intersections cach\'ees \dots   
\bs
{}
\vfill
\eject

{\centerline{TABLE DES MATIERES}
\def\tab{\leaders\hbox to 5mm{\hfil.\hfil}\hfill}
\bs
\bs
\noindent
\line{{\bf I. Sous topos dans les espaces topologiques\tab 1}}

\line{1. Reconstitution d'un espace \`a partir de ses ouverts\tab 1}

\line{2. Notion de sous-lieux\tab 2}
\line{3. D\'efinition de sous-lieux\tab 5}
\line{4. R\'eunions et intersections\tab 5}
\line{5. Images directes \tab 6}
\line{6. Images r\'eciproques\tab 9}
\line{7. Sous-lieux ouverts\tab 9}
\line{8. Sous-lieux ferm\'es\tab 12}
\line{9. Intersections des sous-lieux avec les ouverts et les ferm\'es\tab 14}
\line{10. Int\'erieur, ext\'erieur, adh\'erence, fronti\`ere\tab 15}
\line{11. Lieu g\'en\'erique\tab 16}  
\bs
\noindent
{\bf II. Image r\'eciproque d'une r\'eunion\tab 18}
\bs
\noindent
{\bf III. Th\'eorie de la mesure dans les lieux r\'eguliers\tab 23}

\line{1. D\'efinitions\tab 23}
\line{2. Restriction d'une mesure \`a un sous-lieu\tab 24}
\line{3. Sous lieux r\'eduits et additivit\'e de la mesure \tab 27}
\line{4. Support fin d'une mesure\tab 30}
\bs
\noindent
{\bf IV. Zone d'enchev\^etrement\tab 33}

\line{1. Zone d'enchev\^etrement\tab 33}
\line{2. Crit\`eres de compl\'ementarit\'e\tab 34}
\bs
\noindent
{\bf V. Annexe : sous-espaces et sous-lieux d'un espace\tab 35}

\line{1. Sous-lieu associ\'e \`a un sous-espace\tab 35}
\line{2. Intersections de sous-espaces avec des ouverts et des ferm\'es\tab 36}
\line{3. Union de sous-espaces et union de sous-lieux\tab 36}

\eject
\pageno=1
\noindent
{\bf I - SOUS-TOPOS DANS LES ESPACES TOPOLOGIQUES}
\bigskip
\noindent
{\bf 1 - Reconstitution d'un
espace topologique \`a partir de ses ouverts.}

\medskip
\noindent
{\sl Vocabulaire et notation}. On dit qu'un espace $E$ est {\sl irr\'eductible}
s'il est non-vide et s'il n'existe pas dans $E$ deux ouverts non-vides
disjoints. D'autre part, on appelle {\sl point g\'en\'erique} de $E$ tout
point $x$
tel que $\overline {\langle x \rangle} = E$. Enfin, on dit que $E$ est {\sl
sobre} si toute partie ferm\'ee irr\'eductible de $E$ admet un point
g\'en\'erique et
un seul.

\noindent
 (en particulier : 1) tout espace s\'epar\'e est sobre, puisque tout ferm\'e
irr\'eductible se r\'eduit \`a un point,
2) Si $A$ est un anneau commutatif unif\`ere, Spec A est sobre).
\medskip
Tous les espaces consid\'er\'es dans la suite sont suppos\'es sobres. Pour tout
espace $E$, on notera $O(E)$ l'ensemble des ouverts de $E$ ordonn\'e par
inclusion.

\bigskip
\noindent
{\bf Proposition 1}. {\sl Etant donn\'es deux espaces sobres $E$, $F$ et une
application croissante $\varphi : O(E) \fleche O(F)$, les propositions
suivantes sont \'equivalentes :

a) Il existe une unique application continue $f : F \fleche E$ telle que

$f^{-1}(V) =\varphi(V)$ pour tout ouvert $V$ de $E$.

b) $\varphi(\emptyset) = \emptyset$, $\varphi(E) = F$, et $\varphi$ commute aux
intersections finies et aux r\'eunions quelconques.}

\bigskip
\noindent
(En particulier si $F$ se r\'eduit \`a un point, on obtient une description de
l'ensemble des points de $E$ en termes de l'ensemble des ouverts.)

Cette proposition montre qu'on peut plonger la cat\'egorie des espaces (sobres)
et applications continues dans une cat\'egorie plus large de ``lieux" ({\sl
loci en anglais}) \footnote{(*)}{ou encore {\sl locales.}}.

\bigskip
\noindent
{\bf D\'efinition 1}. Un lieu $E$ consiste en la donn\'ee d'un ensemble
ordonn\'e
$(O(E),\subset)$, ([[dont les \'el\'ements sont appel\'es {\sl ouverts}]])  
qui a les propri\'et\'es suivantes :

a) Il existe un plus petit \'el\'ement $O_E$ (ou $\emptyset$) et un plus grand
\'el\'ement $1_E$ (ou $E$ par abus de notation)

b) Toute famille $(V_i)$ d'\'el\'ements de $O(E)$ a une borne sup\'erieure,
 ([[appel\'ee {\sl r\'eunion}]]), et not\'ee 
$\displaystyle{\un_i}V_i$.

c) Tout couple $U, V$ d'\'el\'ements de $O(E)$ a une borne inf\'erieure,
([[appel\'ee {\sl intersection\/} de $U$ et $V$]]), et not\'ee $U\CAP V$.

d) Pour toute famille $(V_i)$ d'\'el\'ements de $O(E)$ et tout $W \in
O(E)$, on a
$$
W \CAP (\un_i V_i) = \un_i (W \CAP V_i)
$$
(En d'autres termes, $O(E)$ est un {\sl treillis distributif complet}).

\bigskip
\noindent
{\bf D\'efinition 2}.\  Etant donn\'es deux lieux $E$, $F$ un morphisme $f : E
\fleche F$ est une application croissante $f^* : O(F) \fleche O(E)$ telle que

a) $f^*(O_F) = O_E$, $f^*(1_F) = 1_E$

b) $f^*$ commute aux intersections finies et aux unions quelconques.
\bigskip
Pour tout $U \in O(E)$, l'ensemble des $V \in O(F)$ tels que $f^*(V) \subset U$
a donc un plus grand \'el\'ement qu'on notera $f_*(U)$.

On a la relation d'adjonction
$$
f^*(V) \subset U \Longleftrightarrow V \subset f_*(U)
$$
Elle implique que $f^*$ est enti\`erement d\'etermin\'ee par $f_*$.

\bigskip
\noindent
{\sl Commentaire}. La structure de lieu est \'equivalente \`a celle de topos (de
Grothen-\break dieck) {\sl engendr\'e par ses ouverts} ; les morphismes de lieux
s'identifient aux morphismes (g\'eom\'etriques) de topos. D'apr\`es la
proposition 1,
on a une \'equivalence naturelle entre la cat\'egorie des espaces sobres et
applications continues, et une sous-cat\'egorie pleine de la cat\'egorie
des lieux.

\vskip 20pt
\noindent
{\bf 2 - Notion de sous-lieu}.

\medskip
\noindent
{\bf Lemme 1 et d\'efinition}.\ {\sl Pour tout morphisme de lieux $f : E \fleche F$,
les propri\'et\'es suivantes sont \'equivalentes

(i) $f^*$ est surjective

(ii) $f_*$ est injective

(iii) $f^*f_* = 1_{O(E)}$\quad [[{\rm Les compositions sont not\'ees par
 juxtaposition.}]]

\noindent
On appellera {\rm plongement} tout morphisme de lieux ayant ces propri\'et\'es.}

\bigskip
\noindent
{\sl Preuve :\/} d\'ecoule imm\'ediatement de la relation d'adjonction entre
$f^*$ et
$f_*$.

\bigskip
\noindent
{\bf Remarque}. Soit $f : E \fleche F$ une application continue entre espaces
sobres. Pour que $f$ soit un plongement (en tant que morphisme de lieux) il
faut et il suffit qu'elle soit un hom\'eomorphisme de $E$ sur un sous-espace de
$F$.

\bigskip
\noindent
{\sl Preuve}. La condition est suffisante par d\'efinition de la topologie
induite. 

Inversement, supposons que $f$ soit un plongement. Soit $x, x'$ deux
points de $E$ tels que $f(x) = f(x')$. Alors $f_*(E - \overline{\langle x
\rangle}) = f_*(E - \overline{\langle x' \rangle})$, d'o\`u $\overline{\langle x
\rangle} = \overline{\langle x' \rangle}$ puisque $f_*$ est injective, et $x =
x'$ puisque $E$ est sobre. Donc l'application $x \mapsto f(x)$ est une
bijection de l'ensemble des points de $E$ sur un ensemble de points de $F$.
Puisque $f^*$ est surjective, c'est un hom\'emorphisme sur le sous-espace
correspondant.

\bigskip
\noindent
{\bf Lemme 2}. {\sl Soient $i : X \fleche E$ un plongement et $f : Y \fleche E$
un morphisme de lieux. Pour que $f$ se factorise par $i$, il faut et il suffit
que $f_* f^*(V) \supset i^{}_* i^*(V)$ pour tout $V \in O(E)$}.
\bigskip
\noindent
{\bf Corollaire :}\ {\sl deux plongements $i : X\fleche E$ et $j : Y\fleche E$ sont
isomorphes si et seulement si $i^{}_ *i^*=j_*j^*$.}
\bigskip
\noindent
{\sl Preuve du lemme 2 :}
\noindent
 {\sl Condition n\'ecessaire}. Pour tout $V \in O(E)$, on a $V
\subset f_*f^*(V)$ par adjonction entre $f_*$ et $f^*$. Soit $g : Y \fleche Y$
tel que $f = i    g$. On a successivement
$$
\eqalign{
V & \subset i^{}_* g_* g^* i^*(V) \cr
i^*(V) & \subset g_* g^* i^*(V)\quad\quad\quad {\hbox{\rm par adjonction}} \cr
i^{}_*i^*(V) & \subset i^{}_* g_*g^*i^*(V) = f_*f^*(V)
}
$$
\noindent
{\sl Condition suffisante}. On a $f^*f_*f^* = f^*$ : en effet, si $V \in O(E)$,
on a $V \subset f_*f^*(V)$, donc $f^*(V) \subset f^* f_* f^*(V)$. D'autre part,
l'inclusion $f_*f^*(V) \subset f_*f^*(V)$ donne par adjonction $f^*f_*f^*(V)
\subset f^*(V)$. Supposons maintenant $i^{}_* i^* \subset f_*f^*$. Alors
$$
f^*(V) = f^*i^{}_*i^*(V)
$$
pour tout $V \in O(E)$ : en effet
$$
V \subset i^{}_* i^*(V) \Longrightarrow f^*(V) \subset f^*i^{}_*i^*(V)
$$
et d'autre part
$$
i^{}_*i^*(V) \subset f_*f^*(V) \Longrightarrow f^* i^{}_* i^*(V) \subset f^*f_*f^*(V)
= f^*(V)
$$
Par cons\'equent, \'etant donn\'es $V, W \in O(E)$, la relation $i^*(V) =
i^*(W)$
implique $f^*(V) = f^*(W)$. Puisque $i^*$ est surjective, il existe donc une
application et une seule $g^* : O(X) \fleche O(Y)$ telle que $g^* i^* = f^*$.
Puisque $f^*$ et $i^*$ commutent aux intersections finies et aux
r\'eunions, il en
est de m\^eme pour $g^*$, cqfd.

\bigskip
\noindent
{\bf Lemme 3}.  {\sl Soient $E$ un lieu et $e : O(E) \fleche O(E)$ une
application croissante. Les propri\'et\'es suivantes sont \'equivalentes

(a) il existe un lieu $X$ et un morphisme $f : X \fleche E$ telque $e = f_*
f^*$.

(a bis) idem avec un  plongement

(b) $e$ est idempotent, pour tout $U, V \in O(E)$,  $e(U) \supset U$, et
 $e(U \CAP V) = e(U) \CAP e(V)$ }.

\bigskip
\noindent
{\sl Preuve}. (a) $\Longrightarrow$ (b). Les deux premi\`eres
propri\'et\'es de b) ont
\'et\'e \'etablies dans la d\'emonstration du lemme pr\'ec\'edent ; d'autre
part $f_*$
commute aux inf quelconques pour tout morphisme de lieux $f$, puisque c'est
l'adjointe \`a droite de $f^*$.

(b) $\Longrightarrow$ (a bis). Soit $\Omega$ l'ensemble des $V \in O(E)$ tels
que $e(V) = V$ ; c'est aussi l'image de l'application $e$. Prouvons que
$\Omega$ est un treillis distributif complet pour l'ordre induit.

1) Si $U, V \in \Omega$, alors $e(U \CAP V) = e(U) \CAP e(V) = U \CAP V$, donc
$U \CAP V \in \Omega$.

2) Soient $(V_i)$ une famille d'\'el\'ements de $\Omega$ et $W$ un
\'el\'ement de
$\Omega$. On a : 
$$\eqalign{\forall i \in I,\quad V_i \subset W & \Longleftrightarrow \un_i V_i \subset W\cr
                                               & \Longleftrightarrow e(\un_i V_i) \subset W\cr
}$$
\noindent
donc $e(\un_i V_i)$ est borne sup\'erieure de la famille $(V_i)$  dans
$\Omega$.

3) Enfin, avec les m\^emes notations, on a
$$
W \CAP e(\un_i V_i) = e(W) \CAP e(\un_i V_i) = e(W \CAP (\un_i V_i)) =
e(\un_i (W \CAP V_i))
$$
d'o\`u la distributivit\'e. 

On peut donc d\'efinir un lieu $X$ par $O(X) =
\Omega$ ;
nous avons aussi d\'emontr\'e que l'application $e : O(E) \fleche \Omega$
qui est
surjective est l'image inverse pour un morphisme de lieux $i : X \fleche E$ qui
est un plongement. Etant donn\'es $V \in O(E)$ et $W \in \Omega$, on a
$$
e(V) \subset W \Longleftrightarrow V \subset W
$$
donc $i^{}_*$ est l'inclusion $\Omega \fleche O(E)$, d'o\`u $i^{}_* i^* = e$.
\bigskip
\eject
\noindent
{\bf 3 -  D\'efinition d'un sous-lieu ; inclusion.}

\smallskip
Un sous-lieu $X$ d'un lieu $E$ consiste en la donn\'ee d'une application
 $$
e^{}_X :
O(E) \fleche O(E)
$$ ayant les propri\'et\'es \'equivalentes  du lemme 3
(projecteur
associ\'e). On ``mat\'erialise" le lieu correspondant par
$$
O(X) =\im (e^{}_X) = \{V \in O(E)\ |\ e^{}_X(V) = V\}
$$
le plongement\ $i^{}_X : X \fleche E $
\'etant donn\'e par\ $i^*_X(V) = e^{}_X(V),\ (i^{}_X)_*(U) = U$.
\medskip
 Le lemme 2 justifie la
d\'efinition de {\sl l'inclusion} entre sous-lieux : 

$$X \subset Y\quad{\rm si}\quad e^{}_X
\supset
e^{}_Y$$
\bigskip
\noindent
{\bf 4 - R\'eunions et intersections.}
\smallskip
\noindent
{\bf Proposition et d\'efinition.}\ 
{\sl Soient $(X_i)_i$ une famille de sous-lieux de $E$ et $(e^{}_i)_i$ les projecteurs
associ\'es. Pour tout $V \in O(E)$, soit $e(V)$ la r\'eunion des $W \in
O(E)$ qui
sont contenus dans tous les $e^{}_i(V)$. Alors

i) $e$ est le projecteur associ\'e \`a un sous-lieu $X$ de $E$

ii) un sous-lieu $Z$ de $E$ contient $X$ si et seulement si il contient tous
les $X_i$ .

$X$ est encore appel\'e la r\'eunion des $X_i$ et not\'e $\un_i X_i$.}

\bigskip
\noindent
{\sl Preuve de i)}. Puisque $V \subset e^{}_i(V)$ pour tout $i$, on a $V \subset
e(V)$. Pour prouver que $e$ est idempotent, on remarque que
$$
e^{}_i(e(V)) \subset e^{}_i(e^{}_i(V)) = e^{}_i(V)
$$
Enfin, \'etant donn\'es $U, V, W \in O(E)$
$$
\eqalign{
W \subset e(U \CAP V) & \Longleftrightarrow W \subset e^{}_i(U \CAP V) \quad
\hbox{\rm pour tout}\ i  \cr
 & \Longleftrightarrow W \subset e^{}_i(U) \CAP e^{}_i(V) \quad
\hbox{\rm pour tout}\ i  \cr
& \Longleftrightarrow W \subset e(U) \CAP e(V)
}$$

\medskip
\noindent
{\sl Preuve de ii)}. Si $X_i \subset Z$ pour tout $i$, alors $e^{}_Z \subset e^{}_i$
pour tout $i$, donc $e^{}_Z \subset e = e^{}_Y$. Inversement, $X_i \subset Y$ pour
tout $i$ puisque $e \subset e^{}_i$.

On appelle $Y$ la r\'eunion des $X_i$ 
\medskip
On a donc aussi une intersection not\'ee $\CAP X_i$ : {\sl c'est la r\'eunion des
sous-lieux contenus dans tous les } $X_i$. Je ne connais pas de formule
g\'en\'erale pour le projecteur associ\'e \`a une intersection (voir
cependant Lemme
10). On d\'emontrera plus loin que
$$
X \CAP (Y \un Z) = (X \CAP Y) \un (X \CAP Z)
$$
quels que soient les sous-lieux $X, Y, Z$ de $E$ (ce n'est pas imm\'ediat). Mais
l'intersection n'est pas distributive par rapport aux r\'eunions infinies
{\footnote{(*)}{Contre exemple : soit $\gamma_{\Bbb R}$ le lieu 
g\'en\'erique de $\Bbb R$ (cf plus loin en 11); comme $\Bbb R$ est sans point isol\'e 
on a $\gamma_{\Bbb R}\cap [\{x\}]=\emptyset$ pour tout $x$ dans $\Bbb R$.
D'o\`u

 $\gamma_{\Bbb R}\cap (\CUP_{x\in\Bbb R}[\{ x\}])=\gamma_{\Bbb R}
\not= \emptyset=\CUP_{x\in\Bbb R}(\gamma_{\Bbb R}\cap [\{x\}])$  .}}.

\bigskip
\noindent
{\bf 5 - Images directes.}

\smallskip
Soit $f : E \fleche F$ un morphisme de lieux. L'application $f_* f^* : O(F)
\fleche O(F)$ est le projecteur associ\'e \`a un sous-lieu $\im (f)$ de $F$.
D'apr\`es
le lemme 2, $\im (f)$ est le plus petit sous-lieu de $F$ par lequel $f$ se
factorise. Pour tout sous-lieu $X$ de $E$, on d\'efinira l'image directe de $X$
par
$$
f(X) =\im (f\, i^{}_X)
$$

[[voir note ci-dessous \footnote{(**)}{on a la formule fort utile
 $$e^{}_{f(X)}=f_*e^{}_Xf^*$$
En effet :\quad\quad\quad $(fi^{}_X)_*(fi^{}_X)^*=f_*(i^{}_X)_*(i^{}_X)^*f^*=f_*e^{}_Xf^*$
   }]]
\bigskip
\noindent
{\bf Lemme 4}. {\sl Les images directes sont transitives : \'etant donn\'es deux
morphismes de lieux $f : E \fleche F, g : F \fleche G$ et un sous-lieu $X$ de
$E$, on a}
$$
(gf)(X) = g(f(X)).
$$

\bigskip
\noindent
{\sl Preuve}. Soit $Y = f(X)$, et consid\'erons le diagramme commutatif
$$
\diagram{
X & \hfl{u}{} & Y \cr
\vfl{i^{}_X}{} & & \vfl{}{i^{}_Y} \cr
E & \hfl{}{f} & F}
$$
Puisque $Y$ est le plus petit sous-lieu de $F$ par lequel $fi^{}_X$ se factorise, on
a $\im (u) = Y$, d'o\`u $u_* u^*(V) = V$ pour tout $V \in O(Y)$. Si $h = g    f
   i^{}_X$, on a donc
$$
h_* h^* = g_* (i^{}_Y)_* u_* u^* (i^{}_Y)^* g^* = g_* (i^{}_Y)_* (i^{}_Y)^* g^*
$$
\def\eq{\Longleftrightarrow}
\bigskip
\noindent
{\bf Lemme 5}.\ {\sl Pour tout morphisme de lieux $f : E \fleche F$ et toute famille
$(X_i)$ de sous-lieux de $E$, on a}
$$
f(\un_i X_i) = \un_i f(X_i)
$$
[[voir une autre preuve dans la note ci-dessous :\footnote{(*)}{
On a $$\forall i,\quad e^{}_{f(X_i)} =f_*e^{}_{X_i}f^*\quad {\rm et}\quad 
        e^{}_{f(\un_i X_i)} =f_*e^{}_{\un_i X_i}f^*   $$
On a par d\'efinition de $e^{}_{\un_i f(X_i)}$ : pour tous $W,W'$ ouverts de $F$, 
$$\eqalign{
W'\subset e^{}_{\un_if(X_i)}(W) & \eq\forall i,\quad W'\subset (e^{}_{f(X_i)})(W)\cr 
                                 & \eq\forall i,\quad  W'\subset (f_* e^{}_{X_i}f^*)(W)\cr
                                 & \eq\forall i,\quad  f^*(W')\subset e^{}_{X_i}(f^*(W))\cr
                                 & \eq\quad  f^*(W')\subset (e^{}_{\un_i X_i}(f^*(W))\cr
                                 & \eq\quad  W'\subset f_*(e^{}_{\un_i X_i})(f^*(W))\cr
                                 & \eq\quad  W'\subset (e^{}_{f(\un_i X_i)})(W)\cr
}$$
d'o\`u :\quad\quad $f(\un_i X_i) = \un_i f(X_i)$.}]]
\bigskip
\noindent
{\sl Preuve}. On d\'efinit la somme disjointe $S$ des lieux $X_i$ par
$$
O(S) = \Pi^{}_i O(X_i)
$$
Les inclusions $u_i : X_i \fleche S$ \'etant donn\'ees par $u^*_i = pr_i$. Les
plongements $X_i \fleche E$ donnent un morphisme $u : S \fleche E$
$$
u^*(V) = (e^{}_{X_i} (V))_i
$$
La r\'eunion des $X_i$ est l'image de $u$ : 

en effet, pour tout $(W_i)_i \in O(S)$
et tout $V \in O(E)$, on a
$$
\eqalign{
V \subset u_*((W_i)_i) & \Longleftrightarrow u^*(V) \subset (W_i)_i \cr
& \Longleftrightarrow (e^{}_{X_i}(V)) \subset (W_i)_i \cr
& \Longleftrightarrow V \subset W_i \quad \hbox{\rm pour tout}\ i
}$$
Donc $u_* u^*(V)$ est le plus grand \'el\'ement de $O(E)$ contenu 
pour tout $i$ dans $e^{}_{X_i}(V)$.

Notons $T$ la somme disjointe des $f(X_i)$ et $v : T \fleche F$ le morphisme
naturel. Consid\'erons le diagramme commutatif
$$
\diagram{
S & \hfl{f'}{} & T \cr
\vfl{u}{} & & \vfl{}{v} \cr
E & \hfl{}{f} & F
}$$
on a
$$
\eqalign{
f(\un_i X_i) & = f(\im u) = \im (f    u)\quad\quad{\rm (lemme\ 4)} \cr
& = \im (v    f') = v (\im f') \quad\quad{\rm (lemme\ 4)} \cr
& = v(T) = \un_i f(X_i)
}
$$

\bigskip
\eject
\noindent
{\bf 6 - Images r\'eciproques.}

\smallskip
Soient $E \mathop{\fleche}\limits^{f} F$ un morphisme de lieux et $Y$ un
sous-lieu de $f$. D'apr\`es le lemme 5, l'ensemble des sous-lieux $X$ de
$E$ tels
que $f(X) \subset Y$ a un plus grand \'el\'ement que nous noterons $f^{-1}(Y)$.
Plus g\'en\'eralement, pour qu'un morphisme de la forme $h : A \fleche E$ se
factorise par $f^{-1}(Y)$, il faut et il suffit que $f    h$ se factorise
par $Y$ : en effet
$$
\im h \subset f^{-1}(Y) \Longleftrightarrow f(\im h) \subset Y
\Longleftrightarrow \im (f    h) \subset Y
$$
Autrement dit, on a le diagramme cart\'esien
$$
\diagram{
f^{-1}(Y) & \hfl{}{} & Y \cr
\vfl{}{} & & \vfl{}{} \cr
E & \hfl{}{} & F
}$$
[[voir note ci-dessous :\footnote{(*)}{
\noindent
 i) Si $i^{}_X$ d\'esigne l'inclusion d'un sous-lieu $X$ dans $E$ et $Z$ un 
autre sous-lieu, on la formule\quad $i^{-1}(Z)=Z\CAP X$\quad
 comme cas particulier d'image inverse.         

\noindent
ii) On v\'erifie aussi sans difficult\'e  la commutation des images 
inverses aux intersections quelconques :
$$
f^{-1}(\CAP_i(X_i))=\CAP_i f^{-1}(X_i)
$$
}]]
\bigskip
\noindent
{\bf 7 - Sous-lieux ouverts.}

\smallskip
Soient $E$ un lieu et $U, H \in O(E)$. Nous noterons\ $e^{}_U(H)$\ 
le plus grand $W \in O(E)$ tel que $W \CAP U \subset H$. On v\'erifie que $e^{}_U$
est le projecteur associ\'e \`a un sous-lieu que nous noterons provisoirement
$[U]$. Les deux lemmes suivants permettrons d'enlever les crochets.

\bigskip
\noindent
{\bf Lemme 6}.\ {\sl a) Pour tout sous-lieu $X$ de $E$ et tout $U \in O(E)$ 
$$
X \subset [U] \Longleftrightarrow e^{}_X(U) = 1_E
$$

b) Etant donn\'es $U, V \in O(E)$ on a
$$
[U \CAP V] = [U] \CAP [V], \quad \quad e^{}_{U \CAP V} = e^{}_U    e^{}_V = e^{}_V   
e^{}_U \quad{\rm et} \quad U \subset V \eq  [U] \subset [V]
$$\par
c) Pour toute famille $(V_i)$ d'\'el\'ements de $O(E)$, on a
$$
\CUP_i [V_i] = [\CUP_i V_i]
$$

d) Pour tout morphisme de lieux $f : E \fleche F$ et tout $V \in O(F)$ on a
$$
f^{-1}([V]) = [f^*(V)]
$$}

\bigskip
\noindent
{\sl Preuve}.

a) Si $X \subset [U]$, alors $1_E = e^{}_U(U) \subset e^{}_X(U)$. Inversement, si
$e^{}_X(U) = 1_E$, alors pour tout $H \in O(E)$
$$
e^{}_U(H) \subset e^{}_X(e^{}_U(H)) = e^{}_X(U \CAP e^{}_U(H)) = e^{}_X(U \CAP H) \subset e^{}_X(H)
$$

b) Si $U \subset V$, alors pour tout $H \in O(E)$
$$
e^{}_V(H) \CAP U \subset e^{}_V(H) \CAP V \subset H
$$
donc $e^{}_V(H) \subset e^{}_U(H)$, d'o\`u $[U] \subset [V]$.

 Inversement, si $[U]
\subset [V]$, alors $1_E = e^{}_V(V) \subset e^{}_U(V)$, donc $U \subset V$.

 Pour tous $W, H \in O(E)$ on a
$$
\eqalign{
W \subset e^{}_{U \CAP V}(H) & \Longleftrightarrow W \CAP U \CAP V\cr
& \Longleftrightarrow W \CAP U \subset e^{}_V(H)\cr
&   \Longleftrightarrow W \subset e^{}_U(e^{}_V(H))}
$$
donc
$e^{}_{U \CAP V} = e^{}_U    e^{}_V$.

 Si $X$ est un sous-lieu de $E$ contenu dans $[U]$
et $[V]$ on a donc pour tout $H \in O(E)$
$$
e^{}_X(H) = e^{}_X e^{}_X(H) \supset e^{}_X e^{}_V(H) \supset e^{}_U e^{}_V(H) = e^{}_{U \CAP V}(H)
$$
donc $[U] \CAP [V] \subset [U \CAP V]$ ; et on a d\'ej\`a prouv\'e l'inclusion
inverse.

c) Posons $X = \un_i [V_i]$ et $U = \CAP_i V_i$. Pour tout $W, H \in O(E)$ on
a
$$
\eqalign{
W \subset e^{}_X(H) & \Longleftrightarrow \forall i,\  W \subset e^{}_{V_i}(H) \cr
& \Longleftrightarrow \forall i,\ W \CAP V_i \subset H \cr
& \Longleftrightarrow  W \CAP U \subset H  \cr
& \Longleftrightarrow W \subset e^{}_U(H)
}$$

d) Posons $U = f^*(V)$. Prouvons que $f([U]) \subset [V]$. Soit $g$ la
compos\'ee
$[U] \fleche E \fleche F$. Pour tout $H \in O(E)$ on a
$$
f^*(e^{}_V(H)) \CAP U = f^*(e^{}_V(H) \CAP V) = f^*(H \CAP V) \subset f^*(H)
$$
donc $f^*(e^{}_V(H)) \subset e^{}_U f^*(H)$ d'o\`u $e^{}_V(H) \subset f_* e^{}_U f^*(H) =
g_*g^*(H)$.

\medskip
Inversement, soit $X$ un sous-lieu de $E$ tel que $f(X) \subset [V]$. Soit $h =
f    i^{}_X$ la restriction de $f$ \`a $X$. On a
$$
1_F = e^{}_V (V) \subset h_* h^*(V)
$$
donc
$$
1_E = f^*(1_F) \subset h^*(V) = e^{}_X (f^*(V))
$$
d'o\`u $X \subset [f^*(V)]$ d'apr\`es (a).

\bigskip
\noindent
{\bf Lemme 7}. {\sl Soient $X$ un sous-lieu de $E$ et $U \in O(E)$. Pour tout
$V \in O(E)$}
$$
V \subset e^{}_X(U) \Longleftrightarrow [V] \CAP X \subset [U].
$$

\bigskip
\noindent
{\sl Preuve}. Soit $i$ le plongement de $X$ dans $E$. On a
$$
\eqalign{
V \subset e^{}_X(U) & \Longleftrightarrow V \subset i^{}_* i^*(U) \cr
& \Longleftrightarrow i^*(V) \subset i^*(U) \cr
& \Longleftrightarrow X \CAP [V] \subset X \CAP [U]\quad\quad
{\hbox{\rm lemme 6, b) et d)}} \cr
& \Longleftrightarrow X \CAP [V] \subset [U].
}$$

Nous pouvons donc identifier $O(E)$ \`a une partie de l'ensemble $SL(E)$ des
sous-lieux de $E$, stable par intersections finies et r\'eunions quelconques, et
nous parlerons simplement d'``ouverts" de $E$. Le lemme 7 permet de faire tous
les ``calculs" dans $SL(E)$.

\eject
\noindent
{\bf 8 - Sous-lieux ferm\'es.}

\smallskip
Etant donn\'e un ouvert $V$ de $E$, nous d\'efinissons le compl\'ementaire
$E-V$ par
$$
e^{}_{E-V}(H) = H \un V \quad \hbox{\rm pour tout ouvert}\ H.
$$

\bigskip
\noindent
{\bf Lemme 8}. {\sl Pour tout sous-lieu $X$ de $E$}
$$
\eqalign{
V \un X = E &\quad \Longleftrightarrow\quad E-V \subset X \cr
V \CAP X = \emptyset &\quad \Longleftrightarrow\quad X \subset E-V
}$$
\noindent
({\bf Corollaire :}\ $(E-U=E-V)\ \Longrightarrow \ U=V$)
\bigskip
\noindent
{\sl Preuve}. a) Pour tout ouvert $H$ de $E$ on a
$$
\eqalign{
e^{}_X(H) \subset e^{}_X(H \un V) & = e^{}_X(H \un V) \CAP e^{}_V(H \un V) \cr
& = e^{}_{X \un V}(H \un V)}
$$
donc si $X \un V = E$, $e^{}_X(H) \subset H \un V = e^{}_{E-V}(H)$.

Inversement, prouvons que $V \un (E-V) = E$. Pour tout ouvert $H$
$$
\eqalign{
e^{}_V(H) \CAP e^{}_{E-V}(H) & = (e^{}_V(H) \CAP V) \un (e^{}_V(H) \CAP H) \cr
& \subset H \un H = H
}$$

\medskip
b) On a $V \CAP (E-V) = \emptyset$. En effet, si $X$ est un sous-lieu contenu
dans $V$ et $E-V$, alors pour tout ouvert $H$ on a
$$
\eqalign{
e^{}_X(H) & = e^{}_X e^{}_X(H) \supset e^{}_V e^{}_{E-V}(H) \cr
& = e^{}_V (H \un V) = E \cr
}
$$
Le reste de la proposition d\'ecoule du lemme suivant appliqu\'e \`a $g = i^{}_X$.

\bigskip
\noindent
{\bf Lemme 9}. {\sl Pour tout morphisme de lieux $g : A \fleche E$, on a}
$$
g^{-1}(E-V) = A - g^{-1}(V).
$$

\bigskip
\noindent
{\sl Preuve}. Soient $h : B \fleche F$ un morphisme de lieux et $W$ un ouvert
de $F$. Pour que $h$ se factorise que $F-W$, il faut et il suffit que
$h^{-1}(W) = \emptyset$. En effet, 

$$
\eqalign{ h\ {\hbox{\rm se factorise par}}\ 
F-W & \Longleftrightarrow h_* h^*(U) \supset e^{}_{F-W}(U) \quad \forall U \in
O(E) \cr
& \Longleftrightarrow h^*(U) \supset h^* e^{}_{F-W}(U) \quad \forall U \in O(F) \cr
& \Longleftrightarrow h^*(U) \supset h^*(U\un W)=h^*(U)\un h^*(W)\quad \forall U \in O(F) \cr
& \Longleftrightarrow h^*(W) = \emptyset
}$$
Par cons\'equent,
$$
\eqalign{ h : B \fleche A\ {\hbox{\rm se factorise par}}
\quad A - g^{-1}(V) & \Longleftrightarrow h^{-1} g^{-1}(V) = \emptyset \cr
& \Longleftrightarrow gh \ \hbox{\rm se factorise par}\quad E-V \cr
& \Longleftrightarrow h \ \hbox{\rm se factorise par}\quad g^{-1}(E-V)
}$$

\bigskip
\noindent
{\bf Lemme 8 bis}.\quad\quad  $X \un (E-V) = E \Longleftrightarrow V \subset X$, 

\quad\quad\quad\quad\quad\quad $X\CAP (E-V) = \emptyset \Longleftrightarrow X \subset V$.

\bigskip
\noindent
{\sl Preuve}. a) $V \subset X$ veut dire : $e^{}_X(H) \subset e^{}_V(H)$,
 ou encore :
$e^{}_X(H) \CAP V \subset H$, pour tout $H \in O(E)$.

$X \un (E-V) = E$ veut dire $e^{}_X(H) \CAP (H \un V) = H$, ce qui revient au
m\^eme.
\medskip
b) Si $X \subset V$, alors $X \CAP (E-V) = \emptyset$ puisque $V \CAP (E-V) =
\emptyset$.

 Si $X \CAP (E-V) = \emptyset$, alors 
$\emptyset = i_X^{-1}(E-V) = X-i_X^{-1}(V)$,

 donc $i_X^{-1}(V) =X$, d'o\`u $X \subset V$.
\bigskip
\noindent
{\bf Remarques :}
\noindent 
a) Pour toute famille $(V_i)_i$ d'ouverts de $E$, on a
$$
E - \un_i V_i = \CAP_i (E-V_i)
$$
En effet, soit $V = \un_i V_i$. Pour tout sous-lieu $X$ de $E$
$$
\eqalign{
X \subset E-V & \Longleftrightarrow e^{}_X(H) \supset H \un V \quad \hbox{\rm
pour tout}\ h \cr
& \Longleftrightarrow e^{}_X(H) \supset H \un V_i \quad \hbox{\rm pour tout}\ h \
\hbox{\rm et tout}\ i \cr
& \Longleftrightarrow X \subset E-V_i \quad \hbox{\rm pour tout}\ i.
}$$
\noindent
b) Quels que soient les ouverts $V, W$ on a
$$
E - (V \CAP W) = (E-V) \un (E-W)
$$
En effet, si $F = E - (V \CAP W)$, on a pour tout ouvert $H$
$$
\eqalign{
e^{}_F(H) = H \un (V \CAP W) & = (H \un V) \CAP (H \un W) \cr
& = e^{}_{E-V}(H) \CAP e^{}_{E-W}(H)
}
$$
\bigskip
Les compl\'ementaires des ouverts ont donc bien les propri\'et\'es
g\'en\'erales des
ferm\'es d'un espace topologique. L'application $U \longmapsto E-U$ est une
bijection
d\'ecroissante entre les ouverts et les ferm\'es, elle transforme les
intersections
finies en r\'eunions finies et les r\'eunions en intersection.
\ms
\noindent
{\bf 9 - Intersections des sous-lieux avec les ouverts et les ferm\'es.}

\smallskip
On n'a pas de formule g\'en\'erale pour le projecteur associ\'e \`a
l'intersection de
deux sous-lieux, cependant
\ms
\noindent
{\bf Lemme 10}. {\sl Soit $X$ un sous-lieu de $E$.}

{\sl Pour tout ouvert $U$ on a :}
\quad\quad\quad\  $
e^{}_{U \CAP X} = e^{}_U e^{}_X
$

et pour tout ferm\'e $F$ : \quad\quad\quad 
$
e^{}_{F \CAP X} = e^{}_X e^{}_F.
$
\bigskip
\noindent
{\sl Preuve}. a) Etant donn\'es des ouverts $V$ et $H$
$$
\eqalign{
V \subset e^{}_{X \CAP U}(H) & \Longleftrightarrow V \CAP X \CAP U \subset H \cr
& \Longleftrightarrow V \CAP U \subset e^{}_X(H) \cr
& \Longleftrightarrow V \subset e^{}_U e^{}_X(H)
}
$$

b) Pour tout ouvert $H$, on a
$$
e^{}_Xe^{}_F(H) \subset e^{}_{X \CAP F}\, e^{}_{X \CAP F}(H) = e^{}_{X \CAP F}(H)
$$
Il suffit donc de v\'erifier que $e = e^{}_X e^{}_F$ est le projecteur associ\'e
\`a un
sous-lieu. Soit $V$ l'ouvert compl\'ementaire de $F$. Pour tout ouvert $H$, on a
$$
\eqalign{
e(H) & = e^{}_X(V \un H) \cr
e^{}_F e(H) & = e^{}_X (V \un H) \un V = e^{}_X(V \un H)
}$$
et donc $ee(H) = e(H)$. Les autres propri\'et\'es sont \'evidentes.
\bigskip
\noindent
[[voir corollaire ci-dessous\footnote{(*)}{{\bf Corollaire :}\ {\sl Soient $V$ et $W$
 deux ouverts de $E$ et $X$ un sous-lieu,
alors on a :} 
$$X\CAP (V\un W)=(X\CAP V)\un(X\CAP W)
$$ 
(La v\'erification est imm\'ediate.)}]] 
\bs
\noindent
{\bf Lemme 11}. {\sl Soient $X, Y, L$ trois sous-lieux de $E$. Si $L$ est
ouvert ou ferm\'e, alors}
$$
L \CAP (X \un Y) = (L \CAP X) \un (L \CAP Y)
$$
\bigskip
\noindent
{\sl Preuve}. Posons $M = L \CAP (X \un Y)$, $N = (L \CAP X) \un (L \CAP Y)$.
Si $L$ est ouvert, on a pour tout ouvert $H$
$$
\eqalign{
e^{}_M(H) = & e^{}_L e^{}_{X \un Y}(H) = e^{}_L(e^{}_X(H) \CAP e^{}_Y(H) \cr
& = e^{}_L e^{}_X (H) \CAP e^{}_L e^{}_Y (H) \cr
& = e^{}_{L \CAP X}(H) \CAP e^{}_{L \CAP Y}(H) = e^{}_N(H)
}
$$
et si $L$ est ferm\'e
$$
\eqalign{
e^{}_M(H) & = e^{}_{X \un Y} e^{}_L(H) = e^{}_X e^{}_L(H) \CAP e^{}_Y e^{}_L(H) \cr
& = e^{}_{X \CAP L}(H) \CAP e^{}_{Y \CAP L}(H) = e^{}_N(H)
}
$$

\bigskip
\noindent
{\bf 10 - Int\'erieur, ext\'erieur, adh\'erence, fronti\`ere.}

\smallskip
On a encore pour un sous-lieu quelconque $X$ de $E$ la tripartition usuelle en
int\'erieur, fronti\`ere et ext\'erieur.

\medskip
\noindent
{\bf D\'efinitions} :

 $\int X$ = plus grand ouvert contenu dans $X$

$\overline X$ = plus petit ferm\'e contenant $X$

$\partial X =\overline X \CAP (E - \int X)$

Ext $X$ = plus grand ouvert de $E$ disjoint de $X$.

\medskip
Commentons la d\'efinition de l'ext\'erieur. On montrera plus loin que
l'ensemble
des sous-lieux $Y$ de $E$ tels que $X \un Y = E$ a un plus petit
\'el\'ement $X^c$.
Mais celui-ci n'est pas n\'ecessairement disjoint de $X$ (i.e. $X$ n'a pas
toujours un ``compl\'ementaire") et son int\'erieur peut \^etre strictement
plus grand
que l'ext\'erieur de $X$.

V\'erifions les relations usuelles entre  $\int X$,  $\ext X$, $\overline X$ et
$\partial X$.

1)  $\overline X = E-\ext X$\quad\quad (imm\'ediat)

2) $\partial X = E - (\int X \un \ext X)$ \  parce que
$$
\eqalign{
\overline X \CAP (E - \int X) & =(E - \ext X) \CAP (E - \int X) \cr
& = E - (\ext X \un \int X)
}$$

3) $\int X\un\partial X=\overline X$ : soit $e$ le projecteur associ\'e \`a 
  $\int X \un (E - \ext X \un \int X))$. On a pour tout
ouvert $H$
$$
\eqalign{
e(H) & = e^{}_{\int X}(H) \CAP (H \un \ext X \un \int X) \cr
& = H \un (e^{}_{\int X}(H) \CAP \ext X) \un (e^{}_{\int X}(H) \CAP \int X)
}$$
Le dernier terme est contenu dans $H$ ; d'autre part, $\ext X \subset
e^{}_{\int X}(H)$\ puisque\ $\ext X \CAP  \int X = \emptyset$. Donc
$$
e(H) = H \un \ext X.
$$

4)  $\ext X \un \partial X = E - \int X$. En effet
$$
\eqalign{
\int X \CAP (\ext X \un \partial X) & = (\int X \CAP \ext X)
 \un (\int X \CAP \partial X) \cr
& = \emptyset
}$$
$$
\eqalign{
\int X \un \ext X \un \partial X & = \ext X \un \overline
X \quad \hbox{\rm d'apr\`es (3)} \cr
& = E
}$$
\par
5)\ $X$\ dense\quad $\eq$\quad $e^{}_X(\emptyset)=\emptyset$.
\bigskip
\noindent
{\bf 11 - Lieu g\'en\'erique, ``l'intersection des compl\'ementaires".}

\smallskip
Pour tout ouvert $H$ de $E$, soit
$$
e^{}_\gamma (H) = \int \overline H = \ext(\ext (H))
$$
L'application $e^{}_\gamma$ est le projecteur associ\'e \`a un sous-lieu $\gamma =
\gamma (E)$\ :

 Il est clair que $H \subset e^{}_\gamma(H)$. Comme $e^{}_\gamma(H)
\subset \overline H$ on a $\overline {e^{}_\gamma(H)} \subset \overline H$ d'o\`u
$e^{}_\gamma e^{}_\gamma = e^{}_\gamma$.

 Prouvons que pour tout couple d'ouverts $H, K$ 
$$
e^{}_\gamma (H \CAP K) = e^{}_\gamma(H) \CAP e^{}_\gamma(K)
$$
\noindent
 Puisque $\overline{H \CAP K} \subset
\overline H \CAP \overline K$, l'inclusion $\supset$ est claire. 

\noindent
Inversement,
soit $V$ un ouvert contenu dans $\overline H \CAP \overline K$. Puisque $V
\subset \overline H$,  alors $H \CAP V$ est dense dans $V$ ; de m\^eme $K \CAP V$ est
dense dans $V$. Donc $H \CAP K \CAP V$ est dense dans $V$, d'o\`u $V \subset
\overline{H \CAP K}$.

\bigskip
\noindent
{\bf Proposition 1 et d\'efinition}.\ {\sl $\gamma(E)$ est le {\sl plus petit
sous-lieu
dense} de $E$. On l'appellera {\sl lieu g\'en\'erique} de $E$.}

\bigskip
\noindent
{\sl Preuve}.  a) $\gamma(E)$ est dense dans $E$. En effet, soit $V$ un ouvert
de $E$ disjoint de $\gamma(E)$. On a $\gamma(E) \subset E-V$, donc
$$
e^{}_\gamma (\emptyset) \supset e^{}_{E-V}(\emptyset) = V
$$
or $e^{}_\gamma(\emptyset) = \emptyset$.

b) Soit $X$ un sous-lieu dense de $E$. Il faut montrer que pour tout ouvert $H$
de $E$, on a $e^{}_X(H) \subset \overline H$. Soit $V = e^{}_X(H) \CAP \ext H$.
On a $V \CAP X \subset H$ et $V \CAP H = \emptyset$, d'o\`u $V \CAP X =
\emptyset$, et $V = \emptyset$.

\medskip
Soit $X$ un sous-lieu contenu dans tous les {\sl ouverts} denses. Soient $H$ un
ouvert de $E$ et $V = H \un \ext H$, qui est dense. On a :

$e^{}_V(H) = {\rm
Ext}\ (\ext H) = \int \overline H = e^{}_\gamma (H)$, donc $e^{}_X(H)
\supset e^{}_\gamma(H)$.

 Par cons\'equent, $\gamma(E)$ est ainsi l'intersection de
tous les ouverts denses.

\bigskip
\noindent
{\bf Lemme 12}. {\sl Pour tout ferm\'e $F$ de $E$, on a}
$$
\gamma(E) \CAP F = \gamma(E) \CAP \int F.
$$

\bigskip
\noindent
{\sl Preuve}. Soit $V$ l'ouvert compl\'ementaire de $F$. D'apr\`es le lemme 9,
$\gamma(E) \CAP F$ est le sous-lieu de $\gamma (E)$ compl\'ementaire de
$\gamma(E) \CAP V$. Soit $W = \int F = \ext V$. Puisque $V \un W$
est dense dans $E$, on a $\gamma(E) \subset V \un W$, donc\footnote{(*)}{
corollaire lemme 10.}
$$
(\gamma(E) \CAP V) \un (\gamma(E) \CAP W) = \gamma (E)
$$
$$
V \CAP W = \emptyset
$$
et d'apr\`es le lemme 8, $\gamma(E) \CAP W$ est le compl\'ementaire de
$\gamma(E)
\CAP V$ dans $\gamma (E)$.

\eject

\noindent
{\centerline{\bf II - IMAGE RECIPROQUE D'UNE REUNION}}

\vskip 20pt
\def\un{\mathop{\cup}\limits}
\def\CAP{\mathop{\cap}\limits}

\noindent
{\bf 1)}\ Il s'agit de d\'emontrer la proposition suivante :

\noindent
{\sl Etant donn\'es un morphisme de lieux $\displaystyle Y\mathop{\fleche}\limits^f X$ et
deux sous-lieux
$A$, $B$ de $X$, on a}
$$
f^{-1}(A \un B) = f^{-1}(A) \cup f^{-1}(B)
$$
{\centerline{------}}
\bigskip\noindent
{\bf Lemme 1}. {\sl Soient $U$ un ouvert de $X$ et $F$ le ferm\'e
compl\'ementaire. Soient
$X' = U\sqcup F$ la somme disjointe de $U$ et $F$,  et $p : X' \fleche X$
 le morphisme canonique.
L'application $A
\longmapsto p^{-1}(A)$ est une bijection entre sous-lieux de $X$ et
sous-lieux de
$X'$.}

\bigskip
\noindent
{\sl Preuve}. 1) Soit $A$ un sous-lieu de $X$. D'apr\`es le lemme 9
appliqu\'e \`a l'inclusion
$A \hookrightarrow X$, on a
$$
A = (A \CAP U) \un (A \CAP F)
$$
or
$$
(A \CAP U) \un (A \CAP F) = p(p^{-1}(A))
$$

2) Soient $B$ un sous-lieu de $U$ et $C$ un sous-lieu de $F$. D'apr\`es le
corollaire du
lemme 10, on a
$$
U \CAP (B \un C) = (U \CAP B) \un (U \CAP C)
$$
donc $ U \CAP (B \un C) = B$
$$
F \CAP (B \un C) = (F \CAP B) \un (F \CAP C)
$$
donc $F \CAP (B \un C) = C$.

Par cons\'equent $p^{-1} p(A') = A'$ pour tout sous-lieu $A'$ de $X'$.

\bigskip
\noindent
{\bf 2) L'alg\`ebre de Boole $b(X)$ engendr\'ee par les ouverts.}

\smallskip
Pour toute famille finie $U_1 ,\ldots, U_n$ d'ouverts de $X$, soit
$At(U_1,\ldots,U_n)$
l'ensemble des sous-lieux non-vides de la forme
$$
H_1 \CAP \ldots \CAP H_n
$$
avec\ $H_i = U_i$ ou $H=X-U_i$.

Soit $b(U_1,\ldots,U_n)$ l'ensemble des r\'eunions (y compris $\emptyset$)
de sous-lieux appartenant \`a
$ At (U_1,\ldots,U_n)$.

\bigskip
\noindent
{\bf Corollaire 1 du lemme 1}.\ {\sl Soient $X'$ la somme disjointe des sous-lieux 
appartenant \`a
$At(U_1,\ldots,U_n)$ et $X' \mathop{\fleche}\limits^p X$ le morphisme naturel.
L'application $A \longmapsto p^{-1}(A)$ est une bijection entre sous-lieux de $X$ et
sous-lieux de $X'$ (et {\rm donc commute aux r\'eunions}).}

\bigskip
\noindent
{\bf Corollaire 2}. {\sl $b(U_1,\ldots,U_n)$ est stable par unions et 
intersections dans
l'ensemble des sous-lieux de $X$, et c'est une alg\`ebre de Boole pour ces
op\'erations}.

\bigskip
\noindent
{\bf Lemme 2}. {\sl Pour tout $H \in b(U_1,\ldots,U_n)$ on a $f^{-1}(H) \in
b(f^{-1}(U_1),
\ldots,f^{-1}(U_n))$, et l'application
$$
\eqalign{
H & \longmapsto f^{-1}(H) \cr
b(U_1,\ldots,U_n) & \fleche b(f^{-1}(U_1),\ldots,f^{-1}(U_n))
}$$
est un morphisme d'alg\`ebre de Boole}.

\bigskip
\noindent
{\sl Preuve}. Si $H$ appartient \`a  $At(U_1,\ldots,U_n)$ alors $f^{-1}(H)$ est vide ou
dans\break
\noindent
 $At(f^{-1}(U_1),\ldots f^{-1}(U_n))$ d'apr\`es le lemme 9.
 Soit $Y'$ la somme
disjointe des
atomes sur $f^{-1}(U_1),\ldots,f^{-1}(U_n)$. On a donc un diagramme cart\'esien
$$
\diagram{
Y' & \hfl{}{} &  X' \cr
\vfl{q}{} & & \vfl{}{p} \cr
Y & \hfl{}{} & X
}$$
d'o\`u le r\'esultat d'apr\`es les deux corollaires pr\'ec\'edents.

\bigskip
\noindent
Soit maintenant $b(X)$ l'ensemble des sous-lieux de $X$ qui sont dans
$b(U_1,\ldots,U_n)$ pour au moins une famille finie d'ouverts $U_1,\ldots,U_n$.

\bigskip
\noindent
{\bf Proposition 1}.  {\sl a) $b(X)$ est stable par intersections et unions
(finies) dans l'ensemble
des sous-lieux de $X$ et c'est une alg\`ebre de Boole.

b) Pour tout $H \in b(X)$, on a $f^{-1}(H) \in b(Y)$, et l'application $H
\longmapsto
f^{-1}(H)$ est un morphisme d'alg\`ebres de Boole $b(X) \fleche b(Y)$.

c) si $A$, $B$ sont deux sous-lieux de $X$ et si $H \in b(X)$, alors $H
\CAP (A \un B)
= (H \CAP A) \un (H \CAP B)$.}

\bigskip
\noindent
{\sl Preuve}. Le point a) r\'esulte de
$$
b(U_1,\ldots,U_n) \un b(V_1,\ldots V_p) \subset b(U_1,\ldots,U_n, V_1,\ldots, V_p)
$$
Le point (b) est le lemme 2. Point (c) : corollaire 1 du lemme 1.

\bigskip
\noindent
{\bf Lemme 3}. {\sl Tout sous-lieu de $X$ est intersection filtrante de
sous-lieux de $b(X)$.}

\bigskip
\noindent
{\sl Preuve}.  Soit $A$ un sous-lieu de $X$. Pour tout ouvert $V$ de $X$, soit 
l'\'el\'ement de $b(X)$ suivant :
$$
D_V = V \un (X - e^{}_A(V))
$$
$D_V$ appartient … $b(X)$ ; on va montrer que
$$
A = \CAP_V D_V
$$
d'o\`u le r\'esultat puisque $b(X)$ est stable par intersections finies.

\medskip
a) $A \subset D_V$ pour tout ouvert $V$ de $X$. On a, pour tout ouvert $W$ de $X$,
$$
\eqalign{
e^{}_{D_V}(W) & = e^{}_V (W) \CAP (W \un e^{}_A(V)) \cr
& = W \un (e^{}_V(W) \CAP e^{}_A(V))
}$$
or $e^{}_A (V) \CAP A \subset V$, donc
$$
e^{}_V(W) \CAP e^{}_A(V) \CAP A \subset e^{}_V(W) \CAP V \subset W
$$
et
$$
e^{}_V(W) \CAP e^{}_A(V) \subset e^{}_A(W)
$$

b) Soit $A'$ l'intersection de $D_V$. Pour tout ouvert $V$, on a
$$
e^{}_{A'}(V) \supset e^{}_{D_V}(V) = V \un (e^{}_V(V) \CAP e^{}_A(V)) = e^{}_A(V)
$$
donc $A' \subset A$.

\bigskip
\noindent
{\bf Lemme 4}. {\sl Soient $A$ un sous-lieu de $X$ et $(B_i)$ une famille
de sous-lieux
de $X$.}

\noindent
{\sl On a\/} \footnote{(*)}{Remarque : si $X$ est
 un lieu,
 l'ensemble des sous-lieux de $X$, not\'e $S(X)$  muni de l'inclusion 
 {\sl n'est pas\/} un lieu en g\'en\'eral : en effet  (cf [I,4])
 l'intersection n'est pas  en g\'en\'eral distributive par rapport aux
 unions quelconques. Le lemme ci-dessus   
montre que $S(X)$ devient un lieu lorsqu'on le munit de la relation d'ordre
 $\supset$. Le contre-exemple de [I,4] montre que ce lieu n'est pas en g\'en\'eral 
bool\'een (contrairement \`a ce qui ce passe en topologie ``ordinaire").

  Les formules d\'ej\`a \'etablies entrainent que l'application qui \`a tout 
ouvert $U$ de $X$ associe $X\setminus [U]$ d\'efinit un morphisme de lieu 
$q_X : S(X)\fleche X$. De plus si $X$ et $Y$ sont deux lieux et $f : X\fleche Y$
un morphisme de lieux l'application ``image inverse" de $S(Y)$ dans $S(X)$
d\'efinit un  morphisme de lieux $F$ rendant commutatif le diagramme   
$$
\diagram{
S(X) & \hfl{F}{} & S(Y) \cr
\vfl{q_X}{} & & \vfl{}{q_Y} \cr
X & \hfl{}{f} & Y}
$$
Autrement dit  $X\mapsto S(X)$  d\'efinit un foncteur de la cat\'egorie des 
lieux dans elle-m\^eme, et la famille  $(q_X : S(X)\fleche X)_X$ un morphisme de ce foncteur
dans le foncteur identit\'e.}  :
$$
A \un (\CAP_i B_i) = \CAP_i  (A \un B_i)
$$

\bigskip
\noindent
{\sl Preuve}. $\subset$ est \'evidente. Supposons $\supset$ fausse. D'apr\`es la
le lemme 3 il existe alors $H \in b(X)$ tel que
$
I= A \un (\CAP_i B_i)$ soit inclus dans $H$ et $J=\CAP_i (A \un B_i)$
 ne soit pas inclus dans $H$
\medskip
Soit $K$ le compl\'ementaire de $H$ dans $b(X)$. On a
$$
I \CAP K = \emptyset \quad J \CAP K \not = \emptyset
$$
d'apr\`es la proposition 1, b appliqu\'ee aux inclusions $I \hookrightarrow X$, $J
\hookrightarrow X$.
Or
$$
\eqalign{
K \CAP (A \un B_i) & = (K \CAP A) \un (K \CAP B_i)\quad {\rm (prop. 1, c)} \cr
& = K \CAP B_i
}$$
d'o\`u \quad $K \CAP J \subset \CAP_i B_i \subset I$,\quad $K = \emptyset$,
\quad $H = X$,\quad $J\subset H$,\quad
contradiction.

\bigskip
\noindent
{\bf Corollaire}.\ {\sl Si $(A_i)_i$ et $(B_j)_j$ sont deux familles de sous-lieux de
$X$, alors}
$$
(\CAP_i A_i) \un (\CAP_j B_j) = \CAP_{i,j} (A_i \un B_j).
$$

\bigskip
\noindent
{\sl Preuve}. Posons $A = \CAP_i A_i$, $B = \CAP_j B_j$.

Le lemme 4 donne
$$
A \un B = \CAP_j  A \un B_j
$$
et aussi $A \un B_j = \CAP_i A_i \un B_j$ d'o\`u le r\'esultat.
\bigskip
{\centerline{--------}} 
\bigskip
\noindent
{\bf 3) Fin de la preuve du r\'esultat principal.}

 Posons :
$$
A = \CAP_i H_i,\quad B = \CAP_j K_j,\quad{\rm avec}\quad  H_i, K_j \in b(X)
$$
On a
$$
\eqalign{
f^{-1}(A) & = \CAP_i f^{-1}(H_i) \cr
f^{-1}(B) & = \CAP_j f^{-1}(K_j) \cr
f^{-1} (A \un B) & = f^{-1} (\CAP_{i,j} (H_i \un K_j))\quad\quad\quad {\rm (cor.\ lemme\ 4)}\cr
& = \CAP_{i,j} f^{-1}(H_i \un K_j) \cr
& = \CAP_{i,j} (f^{-1}(H_i) \un (f^{-1}(K_j)\quad\quad{\rm (prop.1,\ b)}   \cr
& = (\CAP_i f^{-1} (H_i)) \un (\CAP_j f^{-1}(K_j))\quad{\rm (cor.\ lemme\ 4)} \cr
& = f^{-1}(A) \un f^{-1}(B).
}$$

\eject

\centerline{
\bf III - THEORIE DE LA MESURE DANS LES LIEUX REGULIERS}

\vskip 25pt
\noindent
\vskip 30pt
\noindent
{\bf 1. D\'efinitions.} 

{\bf D\'efinition 1.}\ Nous appellerons mesure (born\'ee) sur un lieu $E$ toute
 application\break
$\mu : {\rm Ouv} (E) \fleche \Bbb R^+$ ayant les propri\'et\'es suivantes
$$
\eqalign{
\mu(\emptyset) & = 0 \cr
U \subset V & \Longrightarrow \mu(U) \leq \mu(V) \cr
\mu(U \un V) & = \mu(U) + \mu(V) - \mu(U \CAP V) \cr
\mu (\un_i V_i) & = \SUP_i \mu(V_i) \ \hbox{\rm pour toute famille filtrante
croissante}
}$$
(Si $E$ est un espace \`a {\sl base d\'enombrable}, on d\'emontre que cette notion
est \'equivalente \`a celle de mesure bor\'elienne (positive born\'ee)).

On prolonge la mesure \`a l'ensemble des sous-lieux par
$$
\mu(X) = \INF \{\mu(V)\ |\ V \hbox{\rm voisinage de}\ X\}
$$

\bigskip
\noindent
{\bf Lemme 1}.\ {\sl Pour toute suite croissante $(X_n)_{n\in\Bbb N}$ de sous-lieux de $E$ on a}
$$
\mu( \un_n X_n) = \SUP_n \mu(X_n).
$$

\bigskip
\noindent
{\sl Preuve}. Soit $\varepsilon > 0$, et prenons une suite $(\varepsilon_n)$ \`a
termes $> 0$ telle que $\Sigma \varepsilon_n < \varepsilon$. Pour tout $n$, soit
$V_n$ un voisinage de $X_n$ tel que $\mu(V_n) - \mu(X_n) < \varepsilon_n$.

Posons $W_n = \un^n_{k=0} V_k$, $W = \un_n W_n = \un_n V_n$. On a
$$
\mu(W_n) - \mu(X_n) \leq \sum^n_{k=0} \varepsilon_k
$$
(v\'erification par r\'ecurrence :
$$
\eqalign{
\mu(W_{n+1}) & = \mu(W_n) + \mu(V_{n+1}) - \mu(W_n \CAP V_{n+1}) \cr
& \leq \mu(W_n) + \mu(V_{n+1}) - \mu(X_n) \cr
& \hbox{\rm (puisque} \ W_n \CAP V_{n+1} \  \hbox{\rm est un voisinage de}\ X_n) \cr
& \leq \mu(W_n) - \mu(X_n) + \mu(X_{n+1}) + \varepsilon_{n+1}\hfil.)
}$$
et donc
$$
\mu(W) (= \SUP \mu(W_n)) \leq \varepsilon + \SUP \mu(X_n)
$$
Ceci implique $\mu(\un X_n) \leq \SUP_n \mu(X_n)$ ; l'in\'egalit\'e inverse est
triviale.

\bigskip
\noindent

{\bf D\'efinition 2.} Un lieu est {\sl r\'egulier\/} si pour tout ouvert $U$ de $E$,
 les ouverts $V$ tels que $\overline V \subset U$ recouvrent
 $U$\footnote{(*)}{Si $X$ est un espace r\'egulier alors le lieu associ\'e est 
r\'egulier ; Voir V, corollaire \`a la proposition 5.}.
\medskip
{\bf  Nous supposons d\'esormais que le lieu $E$ est r\'egulier.}
\ms
Le but de cette section est d'\'etablir deux r\'esultats sans \'equivalent
dans la
th\'eorie classique.

1) {\sl Additivit\'e stricte} pour la mesure des sous-lieux
$$
\mu (A \un B) = \mu(A) + \mu(B) - \mu(A \CAP B)
$$

2) {\sl R\'eduction} : Pour tout sous-lieu $A$ de $E$, l'ensemble des
sous-lieux $A'$
de $A$ tel que $\mu(A') =\mu(A)$ a un {\sl plus petit \'el\'ement} $R_\mu(A)$.

\bigskip
\noindent
{\bf 2. Restriction d'une mesure \`a un sous-lieu.}

\bigskip
\noindent
{\bf Lemme 2}. {\sl Tout sous-lieu de $E$ est intersection de ses voisinages}.

\bigskip
\noindent
{\sl Preuve}. Soit $A$ un sous-lieu de $E$. Il suffit de montrer que
$$
e^{}_A(H) = \un_V e^{}_V(H),\quad (V \ \hbox{\rm voisinage de}\ A)
$$
Soit $W$ un ouvert tel que $\overline W \subset K = e^{}_A(H)$. L'ouvert
$$
V = (E - \overline W) \un H
$$
est un voisinage de $A$. En effet
$$
\eqalign{
A \CAP V & = (A \CAP (E - \overline W)) \un (A \CAP H) \cr
& = (A \CAP (E- \overline W)) \un (A \CAP K) \cr
& = A \CAP (K \un (E - \overline W)) = A \CAP E = A.
}$$
Or $W \CAP V = W \CAP H \subset H$ donc
$$
W \subset e^{}_V(H)
$$
Puisque $e^{}_A(H)$ est recouvert par les $W$ tels que $\overline W \subset
e^{}_A(H)$, la proposition est d\'emontr\'ee.

\bigskip
\noindent
{\bf Lemme 3}. {\sl Pour tout ouvert $U$ de $E$, on a}
$$
\mu(U) + \mu(E-U) = \mu(E).
$$

\bigskip
\noindent
{\sl Preuve}. Soit $(V_\alpha)$ la famille des voisinages de $E-U$, et posons
$W_\alpha = {\rm Ext} (V_\alpha)$. Les $W_\alpha$ forment une famille filtrante
croissante dont la r\'eunion est $U$, d'o\`u
$$
\mu(U) = \SUP \mu(W_\alpha)
$$
or $\mu(W_\alpha) + \mu(V_\alpha) \leq \mu(E)$, d'o\`u $\mu(U) + \mu(E-U) \leq \mu(E)$.

\bigskip
\noindent
{\bf Lemme 4}. {\sl Pour tout ouvert $U$ de $E$ et tout sous-lieu $A$ on a}
$$
\mu(A) = \mu(A \CAP U) + \mu(A \CAP (E-U)).
$$

\bigskip
\noindent
{\sl Preuve}. Soit $W$ un voisinage de $A$. La restriction de $\mu$ \`a Ouv$(W)$
est une mesure, donc d'apr\`es le lemme 3
$$
\eqalign{
\mu(W) & = \mu(W \un U) + \mu(W \CAP (E-U)) \cr
& \geq \mu(A \CAP U) + \mu(A \CAP (E-U)
}$$
donc
$$
\mu(A) \geq \mu(A \CAP U)  + \mu(A \CAP (E-U)).
$$

\bigskip
\noindent
{\bf Corollaire 1}. {\sl Etant donn\'es deux ouverts $U$, $V$ de $E$ et un sous-lieu
$A$, on a}
$$
\mu(A \CAP (U \un V)) = \mu(A \CAP U) + \mu(A \CAP V) - \mu(A \CAP U \CAP V)
$$

\bigskip
\noindent
{\sl Preuve}. D'apr\`es le lemme 4 on a
$$
\mu(A \CAP (U \un V)) = \mu(A) - \mu (A \CAP (E-U) \CAP (E-V))
$$
d'autre part
$$
\eqalign{
\mu(A) & = \mu( A \CAP U) + \mu(A \CAP (E-U)) \cr
& = \mu (A \CAP U \CAP V) + \mu(A \CAP U \CAP (E-V)) + \mu(A \CAP (E-U) \CAP V)
\cr
& + \mu (A \CAP (E-U) \CAP (E-V))
}$$
donc
$$
\eqalign{
\mu(A \CAP (U \un V)) & = \mu(A \CAP U \CAP V)) + \mu(A \CAP U \CAP (E-V) +
\mu(A
\CAP (E-U) \CAP V)\cr
& = \mu(A \CAP U) + \mu(A \CAP (E-U) \CAP V) \cr
& = \mu(A \CAP U) + \mu(A \CAP (E-U) \CAP V) + \mu(A \CAP U \CAP V) - \mu(A
\CAP U \CAP V) \cr
& = \mu(A \CAP U) + \mu(A \CAP V) - \mu(A \CAP U\CAP V)
} $$

\bigskip
\noindent
{\bf Lemme 5}. {\sl Pour toute famille filtrante croissante $(V_\alpha)$
d'ouverts de $E$ et tout sous-lieu $A$, on a}
$$
\mu(A \CAP (\un_{\alpha} V_\alpha) = \SUP_{\alpha} \mu(A \CAP V_\alpha)
$$

\bigskip
\noindent
{\bf Corollaire. Restriction d'une mesure.}

 {\sl Soient $A$ un sous-lieu de
$A$ et $i : A \fleche E$ l'inclusion. L'application :}
$$
\eqalign{
V  \longmapsto \mu(i(V)) \cr
{\rm Ouv}\, (A)  \fleche \Bbb R^+
}$$
est une mesure sur $A$.
\medskip
\noindent
{\sl Preuve du lemme 5.}
\bigskip
\noindent
$\geq\quad$ : parce que $A \CAP (\un V_\alpha) \supset A \CAP V_\alpha$

\noindent
$\leq \quad$ : soit $\varepsilon > 0$ et $W$ un voisinage de $A$ tel que $\mu(W) - \mu(A)
< \varepsilon$.

\noindent
Pour tout $\alpha$ on a
$$
\eqalign{
\mu(A \CAP V_\alpha) + \mu(A \CAP (E-V_\alpha)) & = \mu(A) \quad\quad
                                                    {\hbox{\rm (lemme 4)}} \cr
 \mu(W \CAP V_\alpha) + \mu(W \CAP (E-V_\alpha) & = \mu(W)                \cr 
}
$$
donc $\mu(W \CAP V_\alpha) \leq \mu(A \CAP V_\alpha) + \varepsilon$
$$
\eqalign{
\mu(A \CAP U V_\alpha) & \leq \mu(W \CAP (U V_\alpha)) = \SUP \mu(W \CAP
V_\alpha)\cr
 & \leq \SUP \mu(A \CAP V_\alpha) + \varepsilon}
$$

\eject
\noindent
{\bf 3. Sous-lieux r\'eduits et additivit\'e de la mesure.}

\noindent
{\bf Proposition 1 et d\'efinition}.\ {\sl Pour tout sous-lieu $A$ de $E$, l'ensemble
des sous-lieux $A' \subset A$ tel que $\mu(A') = \mu(A)$ admet un {\sl plus
petit \'el\'ement} $R_\mu(A)$ (``$\mu$-r\'eduction"). On dira que $A$ est
$\mu$-r\'eduit
si $A = R_\mu (A)$.}
\medskip
Quitte \`a restreindre la mesure, on peut supposer $A = E$.

\bigskip
\noindent
{\sl Preuve}.  Pour tout ouvert $U$ de $E$, soit $e^{}_\mu(U)$ la r\'eunion des
ouverts $V \supset U$ tels que $\mu(V) = \mu(U)$. Si $V$ et $V'$ sont deux tels
ouverts, alors
$$
\eqalign{
\mu(V \un V') & = \mu(V) + \mu(V') - \mu(V \CAP V') \cr
& =  \mu(U) + \mu(U) - \mu(U) = \mu(V)
}$$
donc d'apr\`es l'axiome des r\'eunions filtrantes
$$
\mu(e^{}_\mu(U)) = \mu(U)
$$
Nous prouvons maintenant que $e^{}_\mu$ est le projecteur associ\'e \`a un sous-lieu.

\medskip
1) $U \subset V \Longrightarrow e^{}_\mu(U) \subset e^{}_\mu(V)$. En effet
$$
\mu(V \un e^{}_\mu(U)) = \mu(V) + \mu(U) - \mu(V \CAP e^{}_\mu(U)) = \mu(V).
$$

\medskip
2) $e^{}_\mu(e^{}_\mu(U)) = e^{}_\mu(U)$ \  \'evident

\medskip
3) $e^{}_\mu(U \CAP V) = e^{}_\mu(U) \CAP e^{}_\mu(V)$, \hbox{\rm en effet}
$$
\eqalign{
\mu(e^{}_\mu(U) \CAP e^{}_\mu(V)) & = \mu(U) +\mu(V) - \mu(e^{}_\mu(U) \un e^{}_\mu(V)) \cr
& \leq \mu(U) +\mu(V) - \mu(U \un V) = \mu(U \CAP V)
}$$
d'o\`u $e^{}_\mu(U) \CAP e^{}_\mu(V) \subset e^{}_\mu(U \CAP V)$ ; l'autre inclusion est
triviale.

\bigskip
Soit $R$ le sous-lieu de $E$ d\'efini par $e^{}_\mu$. On va montrer que les
voisinages de $R$ sont exactement les ouverts $V$ tels que $\mu(V) = \mu(E)$
d'o\`u en particulier $\mu(R) = \mu(E)$.

\medskip
(a) Soit $V$ un voisinage de $R$. Pour tout ouvert $W$, on a
$$
{\rm Int} (W \un (E-V)) = e^{}_V (W) \subset e^{}_\mu(W).
$$
Si $W = V$, cela donne $E = e^{}_\mu(V)$, donc $\mu(V) = \mu(E)$.

\medskip
(b) Soit $V$ un ouvert de $E$ tel que $\mu(V) = \mu(E)$. Pour tout ouvert $H$
de $E$, on a
$$
\eqalign{
\mu(V \CAP H) & = \mu(V) + \mu(H) - \mu(V \un H) = \mu(E) +\mu(H) - \mu(E) \cr
& = \mu(H)
}$$
Si $ H = e^{}_V(W)$, on a $V \CAP H = W \CAP H$, d'o\`u
$$
\mu(e^{}_V(W)) = \mu(W) \quad {\rm et}\quad e^{}_V(W) \subset e^{}_\mu(W).
$$
Soit maintenant $X$ un sous-lieu de $E$ tel que $\mu(X) = \mu(E)$. Pour tout
voisinage $V$ de $X$, on a $\mu(V) = \mu(E)$, donc $V$ est un voisinage de
$R$, d'o\`u $R \subset X$ d'apr\`es le lemme 2.

\bigskip
\noindent
{\bf Lemme 6}. {\sl Pour toute suite d\'ecroissante $(V_n)$ d'ouverts de $E$, on
a}
 $$
\mu(\CAP V_n) = \INF \mu(V_n).
$$

\bigskip
\noindent
{\sl Preuve}. Posons $I = \CAP_n V_n$, $F_n = E-V_n$, $G = \un_n F_n$. On a
$G \un I = E$. En effet
$$
\eqalign{
G \un (\CAP_n V_n) & = \CAP_n (G \un V_n)\quad\quad\quad \hbox{\rm(voir chap. II, lemme 4)}\cr
&  \supset \CAP_n (F_n \un V_n) = E.
}$$
Par cons\'equent $\mu(I) \geq \mu(E) - \mu(G)$.

\noindent
Or $\mu(G) = \SUP \mu(F_n)$ donc
$$
\eqalign{
\mu(I) & \geq \INF (\mu(E) - \mu(F_n)) \cr
& = \INF \mu(V_n) \quad\quad\quad \hbox{\rm(Lemme 3)}
}$$
d'o\`u l'\'egalit\'e puisque $I \subset V_n$ pour tout $n$

\bigskip
\noindent
{\bf Lemme 7 (et principal)}. {\sl Pour toute famille filtrante d\'ecroissante $(A_i)$ de
sous-lieux de $E$, on a}
$$
\mu(\CAP_i A_i) = \INF_i \mu(A_i).
$$

\bigskip
\noindent
{\sl Preuve}. Soit $(V_\alpha)$ la famille filtrante des voisinages des $A_i$.
On a
$$
\eqalign{
\CAP_{\alpha} V_\alpha & = \CAP_i A_i \cr
\INF \mu(V_\alpha) & = \INF \mu(A_i)
}$$
Prenons une suite croissante $(\alpha_n)$ telle que
$$
\lim_{n \fleche \infty} \mu(V_{\alpha_n}) = \INF_\alpha \mu(V_\alpha) \
(\hbox{\rm not\'e}\ \lambda).
$$
Soit $I = \CAP V_n$. D'apr\`es le lemme 6,
$$
\mu(I) = \lambda
$$
On va montrer que $R_\mu(I) \subset V_\alpha$ pour tout $\alpha$, ce qui
ach\`evera la d\'emonstration. Il suffit d'\'etablir
$$
\mu(V_{\beta} \CAP I) = \mu(I) \quad \hbox{\rm pour tout}\ \beta
$$
Or $\mu(V_\beta \CAP I) = \INF_n \mu(V_\beta \CAP V_{\alpha_n})$.

Soit $\gamma\geq \beta, \alpha_n$. On a
$$
\mu(V_\beta \CAP V_{\alpha_n})\geq \mu(V_\gamma) \geq\lambda, \quad\quad\quad {\rm cqfd}.
$$

\bigskip
\noindent
{\bf Th\'eor\`eme 1}. {\sl Quels que soient les sous-lieux $A$, $B$
 de $E$, on a}
$$
\mu(A \un B) = \mu(A) + \mu(B) - \mu(A \CAP B)
$$

\bigskip
\noindent
{\sl Preuve}. Soient $(V_\alpha)$ la famille des voisinages de $A$ et
$(W_\beta)$ celle des voisinages de $B$. On a
$$
\mu(A \un B) = \INF \mu(V_\alpha \un W_\beta)
$$
et d'apr\`es le lemme 7 et le lemme 2
$$
\mu(A \CAP B) = \INF \mu(V_\alpha \CAP W_\beta)
$$
donc
$$
\eqalign{
\mu(A\un B) +\mu(A\CAP B) & = \INF (\mu(V_\alpha \un W_\beta) +\mu(V_\alpha)
\CAP W_\beta)) \cr
& = \INF (\mu(V_\alpha) +\mu (W_\beta) = \mu(A) + \mu(B)\quad\quad\ {\rm cqfd}
}$$

\bigskip
\noindent
{\bf Corollaire}. {\sl Supposons que $E$ soit un espace topologique. Soit $A$ un
sous-espace de $E$ et $B$ le sous-espace compl\'ementaire. Si le sous-topos $A
\CAP B$ est vide (ou plus g\'en\'eralement de mesure nulle) alors $A$ est
mesurable.}

\vskip 30pt
\noindent
{\bf 4. Support fin d'une mesure.}

\smallskip
\noindent
{\bf Th\'eor\`eme 2 - 1ère version}{\footnote{(*)}{Il n'y a pas de 2ème version}. 
{\sl L'ensemble des sous-lieux
$\mu$-r\'eduits de
$E$, ordonn\'e par inclusion, est une {\sl alg\`ebre de Boole compl\`ete}. Elle
constitue donc l'ensemble des ouverts d'un lieu $B(E,\mu)$. L'application $V
\longmapsto R_\mu(V)$ d\'efinit un morphisme $B(E,\mu)
\mathop{\fleche}\limits^\varphi E$. La mesure des sous-lieux d\'efinit une
mesure
sur $B(E,\mu)$.}

La preuve r\'esulte des lemmes suivants.

\bigskip
\noindent
{\bf Lemme 8}. {\sl Toute r\'eunion de sous-lieux $\mu$-r\'eduits de $E$ est
$\mu$-r\'eduite}.

\bigskip
\noindent
{\sl Preuve}.  Soit $(A_i)$ une famille de sous-lieux $\mu$-r\'eduit et $A =
\un_i A_i$. Soit $A' \subset A$ tel que $\mu(A') = \mu(A)$. Pour tout indice
$i$ on a
$$
\eqalign{
\mu(A_i \CAP A') & = \mu(A_i) + \mu (A') - \mu(A_i \un A') \cr
& = \mu(A_i) + \mu(A) - \mu(A) = \mu(A_i).}
$$
donc $A_i \CAP A' = A_i$ puisque $A_i$ est $\mu$-r\'eduit, d'o\`u $A_i \subset A$
et $A' = A$.

\bigskip
\noindent
{\bf Lemme 9}. {\sl De toute famille de sous-lieux r\'eduits de $E$ on peut
extraire une famille d\'enombrable qui a la m\^eme r\'eunion}.

\bigskip
\noindent
{\sl Preuve}. Quitte \`a ajouter les r\'eunions finies, on peut supposer qu'il
s'agit d'une famille {\sl filtrante} croissante $(A_i)$. Soit $(i_n)$ une suite
croissante d'indices telle que
$$
\lim_n \mu(A_{i_n}) = \SUP_i \mu(A_i) \quad (\hbox{\rm not\'e}\ \lambda).
$$
Posons $B =\un_n A_{i_n}$. Je dis que $A_j \subset B$ pour tout $j$. Il suffit
de montrer que $\mu(A_j \CAP B) \geq \mu(A_j)$. Or
$$
\eqalign{
\mu(A_j \CAP B) &\geq \SUP \ \mu(A_j \CAP A_{i_n}) \cr
\mu(A_j \CAP A_{i_n}) & = \mu(A_j) + \mu(A_{i_n}) - \mu(A_j \un A_{i_n}) \cr
& \geq \mu(A_j) +\mu(A_{i_n}) - \lambda
}$$
et $\dis\lim_{n\fleche\infty} \mu(A_{i_n}) - \lambda = 0$.

\bigskip
\noindent
{\bf Corollaire}. {\sl Pour toute famille filtrante croissante  $(A_i)$ de
sous-lieux $\mu$-r\'eduits de $E$ on a}
$$
\mu( \un_i A_i) = \SUP \mu(A_i)
$$
(r\'esulte du lemme pr\'ec\'edent en appliquant le lemme 1).

\bigskip
\noindent
{\bf Lemme 10}. {\sl Pour tout sous-lieu $A$ de $E$, il existe un sous-lieu $B$
tel que}
$$
\eqalign{
\mu(A \CAP B) = 0 \cr
\mu(A \un B) & = \mu(E)
}$$

\bigskip
\noindent
{\sl Preuve}. Soit $(V_n)$ une suite d\'ecroissante de voisinage de $A$
telle que
$$
\mu(A) = \lim \mu(V_n)
$$
et posons
$$
\eqalign{
B_n & = E - V_n \cr
B & = \un_n B_n
}$$
On a $\mu(B) = \SUP_n \mu (B_n) = \SUP_n \mu(E) - \mu(V_n) = \mu(E) - \mu(A)$
donc
$$
\mu(B) + \mu(A) = \mu(E)
$$
Et aussi
$$
\eqalign{
\mu(A \un B)   = \SUP_n \mu(A \un B_n) &  = \mu(A) + \SUP_n \mu(B_n) \cr
& = \mu(A) + \mu(B)
}$$
D'o\`u $\mu(A \CAP B) = 0$ d'apr\`es le th\'eor\`eme 1.
\medskip
\noindent
Etant donn\'es deux sous-lieux $\mu$-r\'eduits $A$, $B$ de $E$, posons
$$
A \CAP B  = R_\mu(A \CAP B) 
 = \hbox{\rm le plus grand sous-lieu $\mu$-r\'eduit contenu dans}\ A \
{\rm et}\ B
$$

\bigskip
\noindent
{\bf Lemme 11}. {\sl Pour tout sous-lieu $\mu$-r\'eduit $A$ et toute famille
filtrante croissante $(B_i)$ de sous-lieux $\mu$-r\'eduits, on a}
$$
A \CAP (\un_i B_i) = \un_i (A \CAP B_i).
$$

\bigskip
\noindent
{\sl Preuve}. D'apr\`es le lemme 9, il existe une famille d\'enombrable $(i_n)$
telle que
$$
\eqalign{
\un_n B_{i_n} & = \un_i B_i \quad \hbox{\rm not\'e} \ B \cr
\un_n (A \CAP B_{i_n}) & = \un_i (A \CAP B_i)  \quad \hbox{\rm not\'e} \ C
}$$
On a alors
$$
\eqalign{
\mu(A \un B) & = \lim_n \mu(A \un B_{i_n}) \cr
& = \lim_n (\mu(A) + \mu(B_{i_n}) - \mu(A \CAP B_{i_n}) \cr
& = \mu(A) + \mu(B) - \mu(C)
}
$$
donc $\mu(C) = \mu(A \CAP B)$, d'o\`u $C = A \CAP B$.
\medskip
\noindent
La preuve du th\'eor\`eme 2 d\'ecoule maintenant de la propri\'et\'e de
distributivit\'e
$$
A \CAP (B \un C) = (A \CAP B) \un (A \CAP C)
$$
valable quel que soient les sous-lieux $A$, $B$, $C$ de $E$ .
\vfill
\eject

\noindent
{\centerline{\bf IV - ZONES D'ENCHEVETREMENT}}

\bigskip
\noindent
{\bf 1. Zones d'enchev\^etrement.}
\noindent
{\bf Proposition 1.}\ {\sl a) Tout sous-lieu d'un lieu bool\'een est ouvert.[[voir preuve
ci-dessous\footnote{(*)}{ Un lieu bool\'een \'etant r\'egulier, tout sous-lieu
est l'intersection de ses voisinages ouverts donc aussi ferm\'es donc tout
sous-lieu est ferm\'e donc ouvert\dots.
}.]]}
\ms
\noindent
b) Tout sous-lieu de $E$ [[lieu quelconque]] est r\'eunion de ses sous-lieux bool\'eens.
\ms
\noindent 
c)Tout sous-lieu bool\'een $B$ est \'egal au lieu g\'en\'erique de son
adh\'erence.
\ms
\noindent
{\sl Preuve de b}. Soient $A$ un sous-lieu de $E$ et $B$ la r\'eunion des 
sous-lieux
bool\'eens de
$A$. Si $B \not = A$, il existerait un sous-lieu $C$ tel que $B \cap C =
\emptyset$
et $A \cap C \not = \emptyset$ (Preuve : soit $V$ un ouvert tel que $e_B(V)
\not =
e_A(V)$. On prend $C = e_B(V) \cap (E -e_A(V))$. On a $\gamma(A \cap C) \cap B =
\emptyset$, or $\gamma(A \cap C) \subset B$ par d\'efinition.
\ms
\noindent
{\sl Preuve de c.} Soit $F$ l'adh\'erence de $B$. Puisque $B$ est dense dans
$F$, on a
$\gamma(F) \subset B$, et donc $\gamma(F)$ est un ouvert de $B$ ; son
compl\'ementaire
dans $B$ est intersection de $B$ et d'un ouvert $V$ de $F$. On a
$$
\emptyset = V \cap B \cap \gamma(F) = V \cap \gamma(F)
$$
donc $V = \emptyset$, d'o\`u $\gamma(F) = B$.
\bs
\noindent
{\bf D\'efinition.} Soient $A$ et $B$ deux sous-lieux de $E$. Appelons {\sl zone d'enchev\^etrement de
$A$ et
$B$} tout ferm\'e $F$ tel que $A \cap F$ et $B \cap F$ soient denses dans $F$.
\bs
\noindent
{\bf Proposition 2.}\ {\sl Il existe une zone d'enchev\^etrement $\varepsilon(A,B)$ qui contient
toutes les
autres.}
\ms
\noindent
{\sl Preuve}. Soit $(F_\alpha)$ la famille des zones d'enchev\^etrement de $A$ et $B$ , et $C$ sa
r\'eunion. L'adh\'erence de $A \cap C$ contient $A \cap F_\alpha$, donc
$F_\alpha$, pour
tout $\alpha$. Même chose pour $B$. Donc $\overline C$ est une 
 zone d'enchev\^etrement, cqfd.
\bs
\noindent
{\bf Proposition 3.}\ $\overline{A \cap B} = \varepsilon(A,B)$.
\ms
\noindent
{\sl Preuve}. Soit $F = \varepsilon(A,B)$. Puisque $A \cap F$ et $B \cap F$ sont
denses dans $B$, ils contiennent $\gamma(F)$, donc $A \cap B$ est dense
dans $F$.
\bs
\noindent
{\bf 2. Crit\`eres de compl\'ementabilit\'e.}
\ms
\noindent
{\bf Proposition et d\'efinition 4.}\ {\sl Soit $X$ un sous-lieu d'un lieu $E$. L'ensemble des
sous-lieux
$Y$ de $E$ tel que $X \cup Y = E$ a un plus petit \'el\'ement $X^c$. Soit $F$
l'adh\'erence de $X \cap X^c$. $X \cap F$ est d'int\'erieur vide
relativement \`a $F$. [[Si $X \cap X^c =\emptyset$ alors on dit alors
que $X$ poss\`ede un compl\'ementaire ($X^c$) ou que $X$ est compl\'ement\'e]].}
\bigskip
\noindent
{\bf Corollaire.}\ {\sl Pour que $X$  ait un compl\'ementaire  il faut et il suffit 
que pour tout ferm\'e $F \not =\emptyset$
tel que $X\cap F$ soit dense dans $F$, on ait
 ${\rm Int}_F(X \cap F) \not = \emptyset$.}
\ms
\noindent
{\sl Preuve de la proposition 4.}  Soit $V$ un ouvert de $E$ tel que $V \cap F
\subset
X \cap F$. Soit $Y$ le compl\'ementaire de $V \cap F$. On a $Y \cup X = E$,
donc $X^c
\subset Y$ et $X^c \cap V \cap F = \emptyset$, or $X^c \cap F$ est dense
dans $F$,
donc $V \cap F = \emptyset$.
\ms
\noindent
{\sl Preuve du corollaire}.
\ms
\noindent
{\sl Condition n\'ecessaire}. Soit $F$ un ferm\'e tel que $X \cap F$ soit
dense dans
$F$. Alors $X^c$ n'est pas dense dans $F$, sinon $\gamma(F) \subset X \cap
X^c$ donc
il existe un ouvert $V$ tel que $V \cap X^c = \emptyset$ et $V \cap F \not =
\emptyset$.
\bigskip
\noindent
{\sl Condition suffisante :\/} imm\'ediat.
\ms
\noindent
 On suppose maintenant que $E$ est un {\sl espace}.
\bigskip
\noindent
{\bf Lemme 1}.\ {\sl Tout sous-lieu compl\'ement\'e est un sous-espace.}
\ms
\noindent
{\bf Lemme 2}.\ {\sl Pour tout sous-espace $X$ de $E$, $X^c$ est le compl\'ementaire
ponctuel de $X$.}
\ms
\noindent
{\bf Proposition 5.} {\sl Pour qu'un sous-lieu de $E$ soit compl\'ement\'e, il
faut et il
suffit que ce soit un sous-espace et qu'il n'ait pas de zone
d'enchev\^etrement avec
son compl\'ementaire ensembliste.}
\ms
\noindent
{\sl Preuve}. Soit $X$ un sous-lieu de $E$. Si $X \cap X^c =\emptyset$, alors tout
point de $E$ est soit dans $X$ soit dans $X^c$. Soit $Y$, (resp $Y'$) le sous-espace
form\'e des
points de $X$ (resp. $X^c$). On a
$$
\eqalign{
Y \cup Y' & = E, \quad Y \subset X \cr
Y \cap Y' & = \emptyset, \quad Y' \subset X^c
}$$
donc $Y = X$ et $Y' = X^c$ ; ceci prouve directement la proposition.

\eject
\def\un{\mathop{\cup}\limits}
\def\CAP{\mathop{\cap}\limits}
\def\eq{\Longleftrightarrow}
\def\fleche{\longrightarrow}
\def\im{\mathop{\rm Im}\limits}
\def\ext{\mathop{\rm Ext}\limits}
\def\int{\mathop{\rm Int}\limits}
\noindent
{\centerline {\bf V - ANNEXE : SOUS-ESPACES ET SOUS-LIEUX D'UN ESPACE}}
\bs
\def\ni{\noindent}
\ni
{\bf 1. Sous-lieu associ\'e \`a un sous-espace.}
\ms
\ni
{\sl Notations.\/} On fixe pour la suite un espace (topologique\dots) $E$.
Le lieu associ\'e
\`a $E$ sera not\'e $[E]$ et de m\^eme, pour tout sous-espace $X$ de $E$,
 le sous-lieu de $[E]$ associ\'e \`a $X$ sera not\'e $[X]$. Si on note $i$
l'inclusion de $X$ dans $E$, le sous-lieu $[X]$ est d\'efini par le projecteur
$i_*i^*$ et donc par la formule :
$$
\forall V, W,{\hbox{\rm ouverts de $E$}}\quad\quad 
W\subset e^{}_{[X]}(V)\eq W\cap X \subset V
$$   
\ni
{\sl Remarque 1 :\/} comparez avec le lemme 7 (ch I).

\ni
{\sl Remarque 2 :\/} la trace sur $X$ d\'efinit une bijection canonique entre 
les ouverts de $[X]$, (i.e. les ouverts de $[E]$ fixes par $e^{}_{[X]}$) et
les ouverts de $X$.
\bs
\ni
{\bf Proposition 1.}\ {\sl Soient $X$ et $Y$ deux sous-espaces de $E$ et
 $U$ un ouvert de $E$ alors :
\ms
\ni
a) $X\subset Y\Longrightarrow [X]\subset [Y]$.
\ms
\ni
b) $X\subset U \eq [X]\subset [U]$. 
\ms
\ni
c) Si $E$ v\'erifie que tout sous-espace est intersection de ses voisinages 
(par exemple si $E$ est s\'epar\'e ou simplement \`a points ferm\'es) alors on
a l'\'equivalence :
$$
X\subset Y\Longleftrightarrow [X]\subset [Y].
$$}
\ni
{\sl Preuve.} Le a) est clair. Pour b) l'implication directe est donn\'ee par a)
; pour la r\'eciproque 
nous avons les \'equivalences :
$$
\eqalign{ [X]\subset [U] & \eq\quad e^{}_{[U]}\subset e^{}_{[X]}\cr
                         & \eq\quad\forall V,W\in O(E),\quad W\subset e^{}_{[U]}(V)
                              \Longrightarrow W\subset e^{}_{[X]}(V)\cr
                         &\eq \quad \forall V,W\in O(E),\quad W\cap U\subset V
                              \Longrightarrow W\cap X\subset V\cr
}$$
En particulier si on prend $W=E$, $V=U$ on a : $X\subset U$.
\ms
\ni
Pour c) nous avons (comme ci-dessus) les \'equivalences :
$$
\eqalign{ [X]\subset [Y] & \eq\quad e^{}_{[Y]}\subset e^{}_{[X]}\cr
                         & \eq\quad\forall V,W\in O(E),\quad W\subset e^{}_{[Y]}(V)
                              \Longrightarrow W\subset e^{}_{[X]}(V)\cr
                         &\eq \quad \forall V,W\in O(E),\quad W\cap Y\subset V
                              \Longrightarrow W\cap X\subset V\cr
}$$
Posant alors $W=E$ on en d\'eduit :
$$ 
[X]\subset [Y] \Longrightarrow\quad \forall V\in O(E),\quad  Y\subset V
                              \Longrightarrow  X\subset V
$$
Et donc si tout sous-espace de $E$ est intersection de ses voisinages on a :

$[X]\subset [Y]\quad \Longrightarrow\quad  X\subset Y$.
\ms
\ni
{\bf Remarque 3 :} si $E$ v\'erifie les hypoth\`eses de b) on a une injection 
(croissante) de
 l'ensemble des sous-espaces de $E$ dans l'ensemble des sous-lieux de $[E]$.
 On 
a en g\'en\'eral plus de sous-lieux que de sous-espaces : par exemple si $E$ est
s\'epar\'e non vide sans point isol\'e alors le lieu g\'en\'erique de $E$ (voir Chap I, par.11)
ne correspond \`a aucun sous-espace car il n'a pas de point et il est dense). 
\bs
\ni
{\bf 2. Intersections avec des ouverts et des ferm\'es.}
\ms
\ni
{\bf Proposition 2.}\ {\sl Soit $U$ un ouvert de $E$, et $X$ un sous-espace.
On a la formule :
$$
[U\cap X]=[U]\cap [X]
$$}

\ni
{\sl Preuve :} On a d'une part (cf I. lemme 10), 
$e^{}_{[U]\cap [X]}=e^{}_{[U]}e^{}_{[X]}$ d'o\`u 
$$
 W\subset e^{}_{[U]\cap[X]}(V)\eq W\cap U\subset e^{}_{[X]}(V)
\eq (W\cap U)\cap X\subset V
$$
que l'on peut encore \'ecrire sous la forme :
 $W\cap (U \cap X)\subset V $ et ceci est \'equivalent \`a 
$W\subset e^{}_{[U\cap X]}(V)$.
\ms
\ni 
Pour les ferm\'es de $E$ on a :
\ms
\ni
{\bf Proposition 3 :}\ {\sl Soit $F=E\setminus U$ un ferm\'e de $E$ alors :

\ni
i) Le sous-lieu $[F]$ est
\'egal \`a $[E]\setminus [U]$. 

\ni
ii) Soit $X$ un sous-espace de $E$, on a $[X\cap F]= [X]\cap [F]$. }
\ms
\ni
{\sl Preuve :} laiss\'ee au lecteur\dots
\bs

\ni
{\bf 3. Union de sous-espaces et union de sous-lieux.}     
\ms
\ni
{\bf Proposition 4.}\ {\sl Soit $(X_i)_{i\in I}$ une famille de sous-espaces de 
$E$, on a la formule :
$$
\un_i[X_i]=[\un_i X_i] 
$$
\ni
(autrement dit : le sous-lieu associ\'e \`a une union de sous-espaces est 
l'union des sous-lieux.)}
\ms
\ni
{\sl Preuve :} par d\'efinition d'une union de sous-lieux on a
$$
\eqalign{
\forall U, V\in O(E)\quad\quad V\subset  e^{}_{\un_i [X_i]}(U)\quad & \eq\quad
 \forall i,\quad V  \subset  e^{}_{[X_i]}(U)\cr
                      & \eq\quad \forall i,\quad V\cap X_i\subset U\cr
                     &\eq\quad \un_i(V\cap X_i)\subset U\cr
                     &\eq\quad V\cap(\un_i X_i) \subset U\cr
                    & \eq\quad V\subset e^{}_{[\un_i X_i]}(U)\cr}
$$
{\bf Remarque 4 :} Si $X$ et $Y$ sont deux sous-espaces on a seulement :
$$
[X\cap Y]\subset [X]\cap [Y]
$$
Par exemple dans $\Bbb R$, les sous-espaces $\Bbb Q$ et $\Bbb R\setminus \Bbb Q$
sont ensemblistement disjoints mais \'etant deux parties denses de $\Bbb R$,
les sous-lieux associ\'es $[\Bbb Q]$ et $[\Bbb R\setminus\Bbb Q]$ sont denses 
dans $[\Bbb R]$ et ont donc
une intersection  dense (car elle contient le lieu g\'en\'erique\dots). 
\bs
\ni
{\bf Proposition 5.}\ {\sl Soit $X$ un sous-espace de $E$ alors on a :
$$
\ext [X]= [\ext X],\quad\overline{[X]}=[\overline{X}],\quad 
[\int X]\subset \int [X],\quad{\rm et}\quad \partial [X]\subset [Fr(X)].
$$
\ni
De plus si $E$ v\'erifie que tout sous-espace est intersection de ses voisinages
alors on  a $\int [X] = [\int X]$ et $\partial [X] =[Fr(X)]$.
\ms
\ni
($\partial [X]= \overline{[X]}\cap(E\setminus \int [X]) $ (voir Chap. I, par. 10) ;
 $Fr(X)=\overline{X}\cap(E\setminus \int X)=\overline{E}\setminus \int E$).} 
\ms
\ni
{\sl Preuve :} 
$\ext X$ est le plus grand ouvert $U$ tel que $U\cap X=\emptyset$ ;
donc la proposition 2 nous dit  pr\'ecis\'ement que cela est \'equivalent \`a
$[U]\cap [X] =\emptyset$, d'o\`u $\ext [X]= [\ext X]$.

\ni
Pour l'adh\'erence on a

\ni
 $\overline{[X]}=[E]\setminus\ext [X]=[E]\setminus [\ext X] =[E\setminus\ext X]
=[\overline{X}]$.      

\ni
L'int\'erieur de $[X]$ est le plus grand ouvert inclus dans $[X]$ ; si $U$ 
est un ouvert de $E$ alors 
$U\subset X\Longrightarrow [U]\subset [X]$ d'o\`u $\int [X]\subset [\int X]$

\ni
L'inclusion  ci-dessus entraine 
$\partial [X]\subset [Fr(X)]$.

\ni
Le cas o\`u $E$ v\'erifie que tout sous-espace est intersection de ses voisinages
est alors clair. 
\bs
\ni
{\bf Corollaire.}\ {\sl Si $E$ est un espace topologique alors :
$$
E\quad {\hbox{\rm  est r\'egulier}}\quad \Longrightarrow \quad [E]\quad {\hbox{\rm  est r\'egulier}}
$$}
\ni
{\bf Application.} L'\'enonc\'e pr\'ec\'edent entraine
que l'on peut appliquer les th\'eor\`emes du chapitre III aux mesures
 bor\'eliennes finies positives sur les espaces r\'eguliers ; par exemple
\`a $\Bbb S^1$ muni de la mesure angulaire sur les ouverts \dots

\end